\documentclass[12pt]{article}
\usepackage[a4paper,top=4cm,bottom=4cm,left=3cm,right=3cm]{geometry}
\usepackage[colorlinks=true,citecolor=blue]{hyperref}
\usepackage{mathptmx, amsmath, amssymb, amsfonts, amsthm, mathptmx, enumerate, color,mathrsfs}
\usepackage[dvipsnames]{xcolor}
\usepackage[mathscr]{eucal}
\usepackage{graphicx}
\graphicspath{{plots/}}
\usepackage{caption} 
\usepackage{subcaption}
\usepackage[percent]{overpic}
\usepackage{pict2e}

\usepackage{multirow}
\usepackage{booktabs}
\usepackage{epstopdf}
\usepackage{multicol}
\usepackage{epstopdf}
\usepackage{amsmath}
\usepackage{xcolor}
\usepackage{bm}
\usepackage[most]{tcolorbox}
\usepackage{tikz}
\usepackage{etoolbox} 
\usepackage{listofitems} 

\usepackage{algorithm}
\usepackage{algpseudocode}

\tikzset{>=latex} 
\colorlet{myred}{red!80!black}
\colorlet{myblue}{blue!80!black}
\colorlet{mygreen}{green!60!black}
\colorlet{mydarkred}{myred!40!black}
\colorlet{mydarkblue}{myblue!40!black}
\colorlet{mydarkgreen}{mygreen!40!black}
\tikzstyle{node}=[very thick,circle,draw=myblue,minimum size=22,inner sep=0.5,outer sep=0.6]
\tikzstyle{connect}=[->,thick,mydarkblue,shorten >=1]
\tikzset{ 
  node 1/.style={node,mydarkgreen,draw=mygreen,fill=mygreen!25},
  node 2/.style={node,mydarkblue,draw=myblue,fill=myblue!20},
  node 3/.style={node,mydarkred,draw=myred,fill=myred!20},
}

\newtheorem{theorem}{Theorem}[section]

\theoremstyle{definition}

\newtheorem{remark}[theorem]{Remark}
\numberwithin{equation}{section}

\title{\large\sc On an adjoint-based numerical approach for time-dependent optimal control problems of biomedical interest}
\author{Zahra Mirzaiyan$^{*}$, Pierfrancesco Siena$^{*}$, Pasquale Claudio Africa$^{*}$,\\ Michele Girfoglio$^{\dag}$, Gianluigi Rozza$^{*}$}
\date{\vspace{-1cm}}

\begin{document}

\maketitle

\begin{center}
\small $^{*}$SISSA, International School for Advanced Studies, Mathematics Area, mathLab, via Bonomea 265, I-34136 Trieste, Italy \\
\small $^{\dag}$University of Palermo, Department of Engineering,  Viale delle Scienze, Ed. 8, 90128 Palermo, Italy 
\end{center}


\begin{abstract}

This work develops a rigorous numerical framework for solving time-dependent Optimal Control Problems (OCPs) governed by partial differential equations, with a particular focus on biomedical applications. The approach deals with adjoint-based Lagrangian methodology, which enables efficient gradient computation and systematic derivation of optimality conditions for both distributed and concentrated control formulations. The proposed framework is first verified using a time-dependent advection-diffusion problem endowed with a manufactured solution to assess accuracy and convergence properties. Subsequently, two representative applications involving drug delivery are investigated: (i) a light-triggered drug delivery system for targeted cancer therapy and (ii) a catheter-based drug delivery system in a patient-specific coronary artery. Numerical experiments not only demonstrate the accuracy of the approach, but also its flexibility and robustness in handling complex geometries, heterogeneous parameters, and realistic boundary conditions, highlighting its potential for the optimal design and control of complex biomedical systems.
\end{abstract}
 

\vskip 1cm

\section{Introduction}

In recent years, significant improvements have been made in the development of computational tools designed to support clinicians in disease diagnosis, treatment, and prevention, as well as in the identification of patient-specific healthcare solutions: 
see, e.g., \cite{viceconti2015big, thomas2021patient,sankaran2012patient,lassila2013reduced,manzoni2012model,fevola2021,d2012variational,Formaggia,quarteroni2016geometric,Romarowski2018}.

Mathematical models adopted in biology and medicine, such as Navier-Stokes equations to simulate blood flow patterns in complex 3D cardiovascular geometries  \cite{sienadata,balzotti2024reduced,siena2025hybrid,ballarin2016,ballarin2017numerical,rathore2025projection,fresca2020deep} or the advection-diffusion equation for modeling drug delivery  \cite{ferreira2023,ferreira2023prepub,Ferreiral2022}, typically could not be analytically addressed and therefore require advanced numerical solution techniques. 
In this context, several model parameters, such as boundary conditions, geometric features, or fluid properties, are often uncertain and require a careful calibration, as they can strongly invalidate the correctness of the predicted patterns.
In some situations, this selection process is conveniently based on patient-specific in vivo measurements, such as the blood velocity and/or pressure data in cardiovascular applications: see, e.g., \cite{ballarin2016,ballarin2017numerical,girfogliononintrusive}. In this case, the model parameters can be properly tuned with the aim of minimizing the mismatch between the numerical solution and the available clinical data and obtaining a predictive computational model. 
Another possible scenario is when, driven by a medical therapy, one aims to obtain a desired distribution related to a certain quantity of interest, for instance, a drug. In such cases, model parameters can be tuned to reproduce the target distribution which, somehow, needs to be assimilated in the modeling framework.

In this work, we propose a numerical approach to support the process selection of some key model parameters based on Optimal Control Problems (OCPs). 
In particular, we consider drug delivery problems modelled by time-dependent advection-diffusion equations. 
The data assimilation is performed by constraining the minimization of a proper objective functional measuring the distance between the numerical solution and the desired output. In this context, the model parameter to be tuned is known as control variable. 

To numerically solve the formulated OCPs, we adopt an adjoint-based Lagrangian approach \cite{gunzburger2002,quarteroni2009,hinze2009,ito2008}. This strategy models the optimal control problem at hand as an unconstrained minimization
problem, whose solution corresponds to the minimum of a properly defined Lagrangian functional.

To validate our methodology, we consider several case studies: an academic benchmark and two applications of biomedical interest. Concerning the academic benchmark, we consider a variant of the test case proposed in \cite{fu2009}. Here we have an analytical solution allowing us to perform a rigorous convergence analysis.  As a first biomedical test case, we consider an idealized light-triggered drug delivery system. Such a technique was demonstrated to have a great potential to improve cancer treatment. Previous studies \cite{Ferreiral2022,ferreira2023prepub,ferreira2023} have extensively investigated the direct problem,  consisting of computing the drug profile for an assigned light pulse. Moreover, the investigation was restricted to a one-dimensional configuration. In our analysis, we revisit the problem from an optimal control perspective where the light intensity is the control variable to be properly tuned in order to obtain the desired distribution of drug concentration, which, in principle, could be derived from a medical therapy. 
We also extend the modelling to a two dimensional setting. 
The last application involves a patient-specific geometry of a coronary artery in which a cylindrical catheter is implanted, delivering a controlled amount of drug into the bloodstream. In \cite{bazilevs2007} a numerical investigation of the drug delivery dynamics is carried out starting from an assigned injection value.  On the other hand, we reinterpret this in an optimal
control framework, where the goal is to properly calibrate the injection value to minimize the discrepancy between the model solution and a target drug distribution. 

The rest of the paper is organized as follows. Sec. \ref{ocp} presents the mathematical formulation of the OCP for each test case. 
In particular, Sec. \ref{sec:2_1}  focuses on the academic benchmark,  while Sec. \ref{sec:2_2} and Sec. \ref{sec:2_3} address the light-driven drug delivery system application and the coronary artery problem, respectively.  
Then Sec. \ref{sec:numerical-method} provides some insights about the employed numerical algorithm and  
Sec. \ref{sec:computational_results} presents the numerical results.
Finally, some concluding remarks are drawn in Sec. \ref{sec:concl}.




\section{Problem statement}
\label{ocp}

In general, an OCP is composed of three elements  \cite{bewley2001,Gad-elHak2003,gunzburger2002}:
\begin{itemize}
    \item[-] an \emph{objective functional} to minimize,
    \item[-] a \emph{control variable} to properly set in order to minimize the objective functional,
    \item[-] the \emph{mathematical model} which represents a set of constraints for the optimization step.
\end{itemize}
 

We deal with three different test cases: an academic benchmark (Sec. \ref{sec:2_1}) and two biomedical applications for drug delivery, the former associated with a light-triggered system for cancer treatment (Sec. \ref{sec:2_2}), the latter related to a cylindrical catheter placed in a patient-specific coronary artery (Sec. \ref{sec:2_3}). In all the cases, the \emph{objective functional} represents the mismatch between computational results and specific target data. On the other hand, the \emph{control variable} and the \emph{mathematical model} change based on the specific problem considered.


\subsection{Academic benchmark}\label{sec:2_1}

We consider a variant of the test case proposed in \cite{fu2009}, which involves an optimization problem for a transient advection-diffusion equation. 

Let us consider a fixed spatial computational domain \( \Omega \subset \mathbb{R}^2 \)
and a time interval of interest $\left(0, T_f\right]$. 
The state equation, subject to appropriate boundary and initial conditions, is given by: 
\begin{equation}
\begin{cases}
\partial_t y + \textbf{V} \cdot \nabla y - \epsilon\ \nabla^2 y  = f + u &  \text{in} \  \Omega \times (0,T_f], \\
y = y_D &  \text{on} \  \partial \Omega \times (0,T_f], \\
y = y_0 &  \text{in} \  \Omega \times \{0\}, \\
\end{cases}
\label{eq:manufactured}
\end{equation}
where $\partial_t$ denotes the time derivative, 
 $y$ is the state variable, $u$ is the control variable, $\textbf{V}$ is a prescribed velocity field, $f$ is an external source term, $ \epsilon$ is the diffusion coefficient, $y_0$ is the initial datum and $y_D$ is the Dirichlet boundary value. The objective functional $J$ consists of two contributions: the kinetic energy $J_u$ of the control $u$, acting as a stabilization term \cite{strazzullo2021model,dedeoptimal2007}, and the distance $J_y$ between the state variable $y$ and the target data $y_d$:
\begin{equation}
  J(y, u) =  J_u+ J_y = \frac{\beta_1}{2} \int_0^{T_f}\int_{{\Omega}} u^2 \ d\Omega\ dt +\frac{\beta_2}{2} \int_0^{T_f} \int_\Omega (y-y_d)^2   \ d\Omega\ dt,
\label{eq:J}\end{equation} where $\beta_1 >  0$ and $\beta_2 > 0$ are proper weights. The OCP then reads: Find $u$ such that the functional \eqref{eq:J} is minimized under the constraint that $y$ satisfies the system \eqref{eq:manufactured}. Note that we are dealing with a distributed OCP since the control $u$ is applied throughout the domain $\Omega$. 

To solve the problem, we apply the adjoint-based Lagrangian method \cite{gunzburger2002,hinze2009,ito2008,quarteroni2009}.  In such approach, the solution corresponds to the minimum of a proper Lagrangian functional defined as follows:
\begin{align}\label{eq:L1}
  \mathcal{L}(y, u, \lambda, m_0, m_D) = J(y,u) -  \int_0^{T_f} \int_\Omega \lambda \ (\partial_t y + \textbf{V} \cdot \nabla y - \epsilon\ \nabla^2 y-f-u)\ d\Omega \ dt  \ \nonumber\\
  - \int_{\Omega} m_{0} \ \left(y(t=0)-y_0\right)\ d\Omega - \int_0^{T_f} \int_{\partial \Omega} m_D \ (y-y_D) \ d\Omega \ dt,
\end{align}
where $\lambda$, $m_0$ and $m_D$  are the so-called adjoint variables related to the state variable, to the initial data, and to the Dirichlet boundary condition, respectively. In practice, 
the optimal solution corresponds to the stationary point of the Lagrangian functional, i.e.,
\begin{align}
\nabla \mathcal{L} (y, u,\lambda, m_0,m_D) = (\mathcal{L}_y, \mathcal{L}_u, \mathcal{L}_\lambda,\mathcal{L}_{m_0},\mathcal{L}_{m_D}) = \mathbf{0}.
\end{align} 
By setting to zero the derivatives of $\mathcal{L}$ with respect to the adjoint variables, i.e., by taking $\mathcal{L}_{\lambda} = 0$, $\mathcal{L}_{m_0} = 0$, and $\mathcal{L}_{m_D} = 0$,  we re-obtain  the state equation, initial and boundary conditions, i.e the system \eqref{eq:manufactured}. 
On the other hand, by setting to zero the derivative of $\mathcal{L}$ with respect to the state variable, i.e., by taking $\mathcal{L}_{y} = 0$, we get the equation for the adjoint variable $\lambda$ as well as its boundary and initial conditions:
\begin{equation}
    \begin{cases}
   \partial_t \lambda +\mathbf{V}\cdot \nabla\lambda+\varepsilon \ \nabla^2 \lambda =\beta_2 (y_d-y)  &  \text{in} \ \Omega \times(T_f,0], \\
   \lambda=0  & \text{on}\ \partial\Omega \times (T_f,0], \\
    \lambda=0  &  \text{in}\ \Omega \times \{T_f\}. 
    \end{cases} \label{adjoint-unsteady-transport}
\end{equation}

\begin{remark}\label{remark:forwardbackward}
We observe that the adjoint system needs to be integrated backward in time. This is a general rule in the context of OCPs: see, e.g., \cite{Gad-elHak2003,gunzburger2002, dedeoptimal2007,hounumerical1999,quarteronioptimal,zakia2021}.
\end{remark}
Finally, by setting to zero the derivative of $\mathcal{L}$ with respect to the control variable, i.e., by taking $\mathcal{L}_u = 0$, we get the optimality condition: 
\begin{equation}\label{optimality-unsteady-transport}
\lambda + \beta_1 u=0 \quad \text{in} \ \Omega \times (0,T_f]. 
\end{equation} 
The OCP is given by eqs. {\eqref{eq:manufactured}, \eqref{adjoint-unsteady-transport} and \eqref{optimality-unsteady-transport}.

\begin{remark}\label{remark:adjoint}

The system {\eqref{eq:manufactured}, \eqref{adjoint-unsteady-transport}, \eqref{optimality-unsteady-transport}} is closed with respect to $y$, $u$ and $\lambda$. However, it should be noted that the condition $\mathcal{L}_y = 0$ would also yield the equations for the other two adjoint variables, i.e., $m_0$ and $m_D$, which actually depend explicitly on $\lambda$, i.e., $m_0$ and $m_D$ could be computed once $\lambda$ is known. Since the computation of $m_0$ and $m_D$ is not essential for the purposes of this work, we do not report such additional equations. 
\end{remark}

\subsection{Application 1: Drug delivery in a light-driven system}\label{sec:2_2}

The first test case of biomedical interest that we consider concerns light-triggered drug delivery. 
This technique is based on light-responsive nanocarriers (e.g., polymeric nanoparticles) that carry the drug to the tumor site and release it by the action of an external light stimulus \cite{tang2019,wang2020,barbeiro2017,ji2019,kalaydina2018}. The goal is to obtain controlled delivery of the drug at the cancer location. This is crucial to minimize systemic side effects as well as to maintain the drug concentration within its optimal therapeutic window. 

The direct problem, consisting of computing the drug profile for an assigned light pulse, has been widely investigated in \cite{Ferreiral2022,ferreira2023prepub,ferreira2023} in a one-dimensional framework. Here, we re-consider that from an inverse perspective, i.e., we try to determine the light intensity value that ensures a desired drug profile. 
We refer to both distributed and concentrated OCPs and, in both configurations, the light intensity is the control variable, but:

\begin{itemize}
    \item for the distributed case (Sec. \ref{light:distributed}), the control is applied within the entire domain;
    \item for the concentrated case (Sec. \ref{light:concentrated}), the control is only applied  along the left boundary of the domain, serving as a localized control.
\end{itemize}

Moreover, note that, in the concentrated control case, this work deals with a two-dimensional extension. Anyway, for the sake of generality, the description of the problem refers to a two dimensional framework also for the distributed case. By following \cite{Ferreiral2022,ferreira2023prepub,ferreira2023} we refer to an idealized domain, whose details will be given in Sec. \ref{sec:computational_results}.

\subsubsection{Distributed case}\label{light:distributed}


Let us consider a fixed spatial computational domain \( \Omega \subset \mathbb{R}^2 \) 
and a time interval of interest $(0, T_f]$.
The state equations are given by \cite{ferreira2023prepub}: 

\begin{equation}
    \begin{cases}
    \partial_t c_f =D_d \nabla^2 c_f+ \gamma c_b u & \quad \text{in} \quad \Omega \times (0,T_f], \\
    c_f = c_{f,r} 
    & \quad \text{on} \quad \Gamma_r \times (0,T_f],\\
  \nabla c_f \cdot \mathbf{n} = c_{fN} & \quad \text{on} \quad \Gamma_u \cup \Gamma_l \cup \Gamma_d \times (0,T_f], \\
    c_f = c_{f,0}
    & \quad \text{in} \quad \Omega \times \{0\}, \\
    \partial_t c_b = -\gamma c_b u &  \quad \text{in} \quad \Omega \times (0,T_f], \\
     c_b = c_{b,0} & \quad \text{in} \quad \Omega \times \{0\}.
    \end{cases} \label{light-unsteady-enhancer}
\end{equation}
Here \(\Gamma_l\), \(\Gamma_u\), \(\Gamma_r\), and \(\Gamma_d\) correspond to the left, upper, right, and lower boundary of the domain $\Omega$, respectively. 
  In \eqref{light-unsteady-enhancer} $c_b$ and $c_f$ are the concentration of bounded and free drug, respectively, $u$ is the light intensity, the term $\gamma c_b u$ indicates the quantity of bounded drug which is converted to free drug, 
$D_d$ is the diffusion coefficient of the free drug and $\mathbf{n}$ is the unit outward normal vector.  Moreover, $c_{f,0}$ and $c_{b,0}$ are the initial data and $c_{f,r}$ and $c_{fN}$ are the Dirichlet/Neumann  boundary values. 

We remark that the light intensity $u$ is the control variable of our problem.
The objective functional $J$ involves the sum of two contributions: $J_u$ is related to the energy of the control (ensuring the stability of the problem) 
and $J_{c_{f,T_f}}$ is related to the distance between the free drug concentration and the desired distribution $c_{f,T_f}$ at the final time of the simulation $t = T_f$:
\begin{eqnarray}\label{eq:func_distr}
  J(c_f, c_b, u)=J_u + 
  J_{c_{f,T_f}}=\frac{\beta_1}{2} \int_0^{T_f} \int_{\Omega} u^2  \,d\Omega \ dt +
  \frac{\beta_3}{2} \int_{\Omega} \left(c_f(t=T_f)-c_{f,T_f}\right)^2 d\Omega,
\end{eqnarray}
where $\beta_1 > 0$ 
and $\beta_3 > 0$ are proper weights. 
 The OCP then reads: Find $u$ such that the functional  \eqref{eq:func_distr} is minimized under the constraint that $c_f$ and $c_b$ satisfy the state equations \eqref{light-unsteady-enhancer}.

Let us define the Lagrangian functional as: 
\begin{align}
  \mathcal{L}(&c_f, c_b, u, \lambda_1, \lambda_2, m_{b,0}, m_{f,0}, m_{fr}, m_N)\nonumber\\
  &=J(c_f, c_b, u)-\int_0^{T_f} \int_\Omega \lambda_1(\partial_t c_f-D_d \nabla^2 c_f- \gamma c_b u) \ d\Omega\ dt\nonumber\\
  &-\int_0^{T_f} \int_\Omega \lambda_2(\partial_t c_b+\gamma c_b u) \ d\Omega \ dt- \int_{\Omega} m_{f,0} \ (c_f(t=0)-c_{f,0})\ d\Omega\nonumber\\
  & - \int_{\Omega} m_{b,0} \ (c_b(t=0)-c_{b,0})\ d\Omega-\int_0^{T_f} \int_{\Gamma_r} m_{fr} \ (c_f-c_{f,r})\ d\Omega \ dt \nonumber \\ 
  &-\int_0^{T_f} \int_{\partial \Omega \setminus \Gamma_r} m_N \ (\nabla c_f \cdot \mathbf{n}-c_{fN})\ d\Omega \ dt. \label{lagrangian-light-distributed-boundary}
\end{align}

In equation \eqref{lagrangian-light-distributed-boundary} 
$\lambda_1$ and $\lambda_2$ are the adjoint variables associated to the state variables $c_f$ and $c_b$, respectively. Moreover, we have the adjoint variables,  $m_{b,0}$ and \( m_{f,0} \), \( m_{fr} \) and \( m_N \), related to the initial data and boundary conditions, respectively.
The optimal solution is the one for which we have 
\begin{eqnarray}
\nabla \mathcal{L}(c_f, c_b, u,  \lambda_1, \lambda_2, m_{b,0}, m_{f,0}, m_{fr}, m_N)= \\(\mathcal{L}_{c_f},\mathcal{L}_{c_b}, \mathcal{L}_u, \mathcal{L}_{  \lambda_1},\mathcal{L}_{  \lambda_2},\mathcal{L}_{  m_{b,0}},\mathcal{L}_{  m_{f,0}},\mathcal{L}_{  m_{fr}},\mathcal{L}_{  m_N})=\textbf{0}.\nonumber
\end{eqnarray}
By taking $\mathcal{L}_{\lambda_1} = \mathcal{L}_{\lambda_2} = 0 = \mathcal{L}_{m_{b,0}} = \mathcal{L}_{m_{f,0}} = \mathcal{L}_{m_{fr}} = \mathcal{L}_{m_{N}} = 0$ we get again the system \eqref{light-unsteady-enhancer}. Moreover, by taking $\mathcal{L}_{c_f} = \mathcal{L}_{c_b} = 0$, we obtain the equations for the adjoint variables $\lambda_1$ and $\lambda_2$:  

\begin{equation}
 \begin{cases}
\partial_t \lambda_1 +D_d \nabla^2 \lambda_1=0  &\quad \text{in} \quad \Omega \times  (T_f,0],
 \\
 \lambda_1=0 & \quad \text{on} \quad \Gamma_r  \times (T_f,0],\\
  \nabla \lambda_1 \cdot\textbf{ n}=0 & \quad \text{on} \quad \partial \Omega -\Gamma_r  \times (T_f,0], \\
 \lambda_1=\beta_3 (c_f-c_{f,T_f}) & \quad \text{in} \quad \Omega \times \{T_f\}, \\
  \partial_t \lambda_2+\gamma \lambda_1 u-\gamma \lambda_2 u =0 \ &  \quad \text{in} \quad \Omega \times (T_f,0],\\
 \lambda_2=0, &  \quad \text{in} \quad \Omega \times \{T_f\}.
 \end{cases} \label{adjointdebc-unsteady-light-distributed-lambda1-lambda2}
\end{equation}

For the equations related to the other adjoint variables 
see Remark \ref{remark:adjoint}. Finally, by enforcing $\mathcal{L}_u = 0$, we obtain the optimality condition:
\begin{equation}\label{eq:optimality}
\beta_1 u +\gamma c_b \lambda_1-\gamma c_b \lambda_2=0 \  \quad \text{in} \quad\Omega \times (0,T_f].
\end{equation}

The OCP is given by eqs. 
\eqref{light-unsteady-enhancer}, \eqref{adjointdebc-unsteady-light-distributed-lambda1-lambda2} and \eqref{eq:optimality}.  


\subsubsection{Concentrated case}\label{light:concentrated}

Let us consider a fixed spatial computational domain \( \Omega \subset \mathbb{R}^2 \) 
and a time interval of interest $(0, T_f]$. The state equations are given by \cite{ferreira2023}: 

\begin{equation}
    \begin{cases}
    \frac{1}{\beta} \partial_t I = D_I  \nabla^2 I - \mu_a I & \quad \text{in} \quad \Omega \times (0,T_f], \\
    I = u & \quad \text{on} \quad \Gamma_l \times (0,T_f], \\
    \nabla I \cdot \mathbf{n} = I_N & \quad \text{on} \quad \partial \Omega - \Gamma_l \times (0,T_f],\\
    I = I_{0} & \quad \text{in} \quad \Omega \times \{0\}, \\
    \partial_t c_f =D_d \nabla^2 c_f+ \gamma c_b I & \quad \text{in} \quad \Omega \times (0,T_f], \\
    c_f = c_{f,r} & \quad \text{on} \quad \Gamma_r \times (0,T_f],\\
    \nabla c_f \cdot \mathbf{n} = c_{fN} & \quad \text{on} \quad \partial \Omega - \Gamma_r \times (0,T_f],\\
        c_f = c_{f,0} & \quad \text{in} \quad \Omega \times \{0\}, \\
     \partial_t c_b = -\gamma c_b I & \quad \text{in} \quad \Omega \times (0,T_f], \\
    c_b = c_{b,0} & \quad \text{in} \quad \Omega \times \{0\}. \\
    \end{cases} \label{light-concentrated-unsteady-enhancer}
\end{equation}

Unlike the distributed case, here we have an equation governing the light intensity $I$ (the first of the system \eqref{light-concentrated-unsteady-enhancer}) and we limit to control its value on the left boundary $\Gamma_l$. 
In addition, $\beta$ is the light speed, $\mu_a$ is the absorption coefficient 
and $D_I$ is the light diffusion coefficient. 
Similarly to the distributed case, the objective functional $J$ involves the sum of two contributions: $J_{u}$ is the energy of the control and 
$J_{c_{f,T_f}}$ is the mismatch between the free drug concentration and the desired distribution $c_{f,T_f}$ at the final time of the simulation $t = T_f$: 

\begin{eqnarray}\label{eq:func_concen}
  J(c_f, c_b, u)=J_{u} 
  + J_{c_{f,T_f}}=\frac{\beta_1}{2} \int_0^{T_f}\int_{\Gamma_l} u^2  \ d\Omega \ dt 
  + \frac{\beta_3}{2} \int_{\Omega} \left(c_f(t=T_f)-c_{f,T_f}\right)^2 d\Omega,
\end{eqnarray}
where $\beta_1 > 0$ 
and $\beta_3 > 0$ are proper weights. 
In this case, the OCP reads: Find $u$ 
such that functional \eqref{eq:func_concen} is minimized under the constraint that $I$, $c_f$ and $c_b$ satisfy the system \eqref{light-concentrated-unsteady-enhancer}. 

We define the Lagrangian functional as: 
\begin{align}
\mathcal{L}(&c_f, c_b, I, u, \lambda_1, \lambda_2, \lambda_3, 
m_{I,0}, m_{f,0}, m_{b,0},  m_{fr}, m_{IN}, m_{fN}) 
\nonumber\\
=&\ J(c_f, c_b, u) 
-\int_0^{T_f} \int_{\Omega} 
\lambda_1 \left(\frac{1}{\beta} \partial_t I - D_I  \nabla^2 I + \mu_a I\right) 
\, d\Omega \, dt 
\nonumber\\
&-\int_0^{T_f} \int_{\Omega} 
\lambda_2 (\partial_t c_f -D_d \nabla^2 c_f- \gamma c_b I) 
\, d\Omega \, dt 
-\int_0^{T_f} \int_{\Omega} 
\lambda_3 ( \partial_t c_b +\gamma c_b I) 
\, d\Omega \, dt 
\nonumber\\
&+ \int_{\Omega} m_{I,0}\,(I(t{=}0)-I_0)\, d\Omega
+ \int_{\Omega} m_{f,0}\,(c_f(t{=}0)-c_{f,0})\, d\Omega+ \int_{\Omega} m_{b,0}\,(c_b(t{=}0)-c_{b,0})\, d\Omega\nonumber\\
&+ \int_0^{T_f} \int_{\Gamma_r} m_{fr}\,(c_f-c_{f,r}) \, d\Omega \, dt+ \int_0^{T_f} \int_{\partial \Omega \setminus \Gamma_l} 
m_{IN}\,(\nabla I\cdot \mathbf{n}-I_N)\, d\Omega \, dt\nonumber\\
&+ \int_0^{T_f} \int_{\partial \Omega \setminus \Gamma_r} 
m_{fN}\,(\nabla c_f\cdot \mathbf{n}-c_{fN})\, d\Omega \, dt.
\label{lagrang-light-concen}
\end{align}

 In equation \eqref{lagrang-light-concen},  $\lambda_1$, $\lambda_2$, and $\lambda_3$ are the adjoint variables associated to the state variables $I$, $c_f$ and $c_b$, respectively. Additionally, we also have the adjoint variables, $m_{I,0}$, $m_{f,0}$ and $m_{b,0}$,  $m_{IN}$, $m_{fr}$, and $m_{fN}$, related to initial data and to boundary conditions. 
 The optimal solution is obtained by enforcing: 
 \begin{align}
\nabla \mathcal{L}(&c_f,c_b, I, u,  \lambda_1, \lambda_2, \lambda_3, m_{I,0}, m_{f,0}, m_{b,0},  m_{fr}, m_{IN}, m_{fN})=&\nonumber\\
&(\mathcal{L}_{c_f}, \mathcal{L}_{c_b},\mathcal{L}_{I},\mathcal{L}_{u},\mathcal{L}_{  \lambda_1},\mathcal{L}_{  \lambda_2},\mathcal{L}_{  \lambda_3},\mathcal{L}_{  m_{I,0}},\mathcal{L}_{  m_{f,0}},\mathcal{L}_{  m_{b,0}},\mathcal{L}_{  m_{fr}},\mathcal{L}_{  m_{IN}},\mathcal{L}_{  m_{fN}})=\textbf{0}.
\end{align}

In particular, by taking $\mathcal{L}_{\lambda_1} = \mathcal{L}_{\lambda_2}=\mathcal{L}_{\lambda_3} = 0$ and $\mathcal{L}_{m_{b,0}} = \mathcal{L}_{m_{f,0}} = \mathcal{L}_{m_{I,0}}=\mathcal{L}_{m_{fr}}=\mathcal{L}_{m_{IN}} = \mathcal{L}_{m_{fN}} = 0$, we get the system \eqref{light-concentrated-unsteady-enhancer}. Additionally, by taking $\mathcal{L}_{I}=\mathcal{L}_{c_f} = \mathcal{L}_{c_b} = 0$, we get the adjoint equations for $\lambda_1$, $\lambda_2$ and $\lambda_3$: 
 




\begin{equation}
    \begin{cases}
   \frac{1}{\beta}\partial_t \lambda_1+D_I \nabla^2 \lambda_1-\mu_a \lambda_1- \gamma c_b \lambda_3+ \gamma c_b \lambda_2=0  & \quad \text{in} \quad \Omega \times (T_f,0], \\
   \nabla\lambda_1 \cdot \mathbf{n}=0  & \quad \text{on} \quad \partial \Omega -\Gamma_l \times [T_f,0),\\
   \lambda_1=0  & \quad \text{on} \quad \Gamma_l\times [T_f,0),\\
    \lambda_1=0  & \quad \text{in} \quad \Omega \times \{T_f\}, \\
 \partial_t \lambda_2 +D_d \nabla^2 \lambda_2=0  & \quad \text{in} \quad \Omega \times [T_f,0),\\
   \lambda_2=0  & \quad \text{on} \quad \Gamma_r\times [T_f,0),\\
 \nabla\lambda_2 \cdot \mathbf{n}=0  & \quad \text{on} \quad \partial \Omega -\Gamma_r \times [T_f,0),\\
  \lambda_2=\beta_3 (c_f-c_{f,0})  & \quad \text{in} \quad \Omega \times \{T_f\}, \\
  \partial_t \lambda_3-\gamma \lambda_3 I+\gamma  \lambda_2 I=0 & \quad \text{in} \quad \Omega \times [T_f,0), \\
    \lambda_3=0, & \quad \text{in} \quad \Omega \times \{T_f\}. \\
    \end{cases} \label{adjoint-unsteady-light-lambda1-lambda2-lambda3}
\end{equation}
For the equations related to the other adjoint variables 
see Remark \ref{remark:adjoint}. Finally, by enforcing $\mathcal{L}_{u} = 0$, we obtain the optimality condition: 
\begin{equation}\label{eq:light-optimality-concentrated}
    \beta_1 u-D_I \ \nabla\lambda_1 \cdot \textbf{n}=0  \quad \text{on} \quad \Gamma_l \times (0,T_f]. \\
\end{equation}
The OCP is given by eqs. \eqref{light-concentrated-unsteady-enhancer}, \eqref{adjoint-unsteady-light-lambda1-lambda2-lambda3} and \eqref{eq:light-optimality-concentrated}.

\subsection{Application 2: Drug delivery in a coronary artery}\label{sec:2_3}
The computational domain $\Omega \subset \mathbb{R}^3$ is given by the patient-specific geometry shown in Figure \ref{fig:geometry}. It represents a portion of the coronary arterial tree, namely the Left Main Coronary Artery (LMCA), which bifurcates into the Left Anterior Descending (LAD) artery and the Left Circumflex (LCx) artery. 
A cylindrical catheter is placed next to the inlet section to deliver a certain quantity of drug into the bloodstream. The boundary of the computational domain consists of: the boundary of the catheter where the drug is released denoted by $\Gamma_{\text{drug}}$, the rest of the boundary of the catheter $\Gamma_{\text{c}}$
 , an inlet denoted by $\Gamma_{\text{inlet}}$ where the blood enters the artery, two outlets denoted by $\Gamma_{\text{outlet}}$ where the blood leaves the arterial tree and finally the artery wall denoted by $\Gamma_{\text{wall}}$.

Such a configuration is very similar to the one reported in \cite{bazilevs2007} 
where the drug delivery was treated as a direct problem. On the other hand, we reinterpret this in an optimal control framework. In particular, the goal is to estimate the quantity of drug that should be delivered from $\Gamma_{\text{drug}}$ in order to minimize the mismatch between the computed drug concentration and a given optimal drug distribution, which could be, in principle, derived from a medical therapy.
From an operative point of view, we refer to a concentrated OCP where the control variable is the Dirichlet value on $\Gamma_{\text{drug}}$. 

The blood is modeled as an incompressible Newtonian fluid. The velocity $\mathbf{V}$ and the pressure $p$ of the blood are derived using the unsteady incompressible Navier-Stokes equations:

\begin{figure}[htb!]
    \centering
    \includegraphics[width=0.8\textwidth]{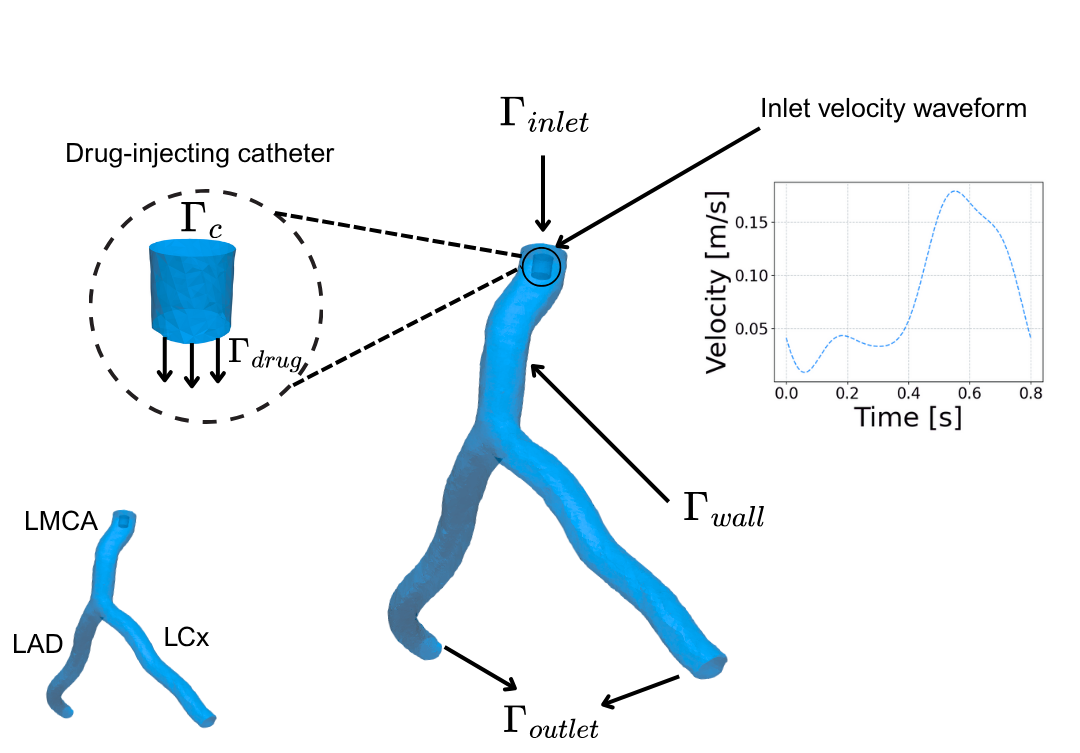}
    \caption{Coronary artery system: sketch of the geometry, including LMCA (Left Main Coronary Artery), LAD (Left Anterior Descending), and LCx (Left Circumflex). }
    \label{fig:geometry}
\end{figure}

\begin{equation}
	\begin{cases} 
    \partial_t \mathbf{V} + (\mathbf{V}  \cdot \nabla) \mathbf{V}
- \nu \nabla^2 \mathbf{V} + \nabla  p  = 0  & \quad \text{in} \quad \Omega \times (0, T_f], \\
	\nabla \cdot \mathbf{V} = 0 & \quad \text{in} \quad \Omega \times (0, T_f],\\
    \mathbf{V} = \mathbf{0} &\quad \text{in} \quad \Omega \times \{0\}, \\
    \mathbf{V} = {V}_i \mathbf{n} &\quad \text{on} \quad \Gamma_{\text{inlet}}\times (0, T_f],\\
    \mathbf{V} = {V}_d \mathbf{n} &\quad \text{on} \quad \Gamma_{\text{drug}}\times (0, T_f],\\
    \mathbf{V}  = \mathbf{0} &\quad \text{on} \quad \Gamma_{\text{wall}} \cup \Gamma_{c}  \times (0, T_f], \\
    (2\mu \nabla \mathbf{V} - p \mathbf{I})\mathbf{n} = \mathbf{0} &\quad \text{on} \quad \Gamma_{\text{outlet}}\times (0, T_f],
	\end{cases} \label{tN-Steady}
\end{equation}
where $\nu$ is the constant kinematic viscosity. 
It is assumed that the drug concentration $y$
does not affect the flow physics, hence we have a one-way coupling and  the drug is modeled as a passive scalar whose governing equation is given by: 


\begin{equation}
    \begin{cases}
    \partial_t y+\mathbf{V} \cdot \nabla y- \varepsilon\ \nabla^2 y= 0& \quad \text{in} \quad \Omega \times (0, T_f],\\
    y= y_0  & \quad \text{in} \quad \Omega \times \{0\},\\
     y= u  & \quad \text{on} \quad \Gamma_{\text{drug}} \times (0, T_f], \\
     y = y_{\text{c}}  & \quad \text{on} \quad \Gamma_{\text{c}} \times (0, T_f], \\
     y= y_{\text{in}}  & \quad \text{on} \quad \Gamma_{\text{inlet}} \times (0, T_f], \\
     \nabla y\cdot \textbf{n} = y_{\text{N}} & \quad \text{on} \quad \Gamma_{\text{outlet}} \cup \Gamma_{\text{wall}}  \times (0, T_f], \\
     \end{cases} \label{eq:general_transport}
\end{equation}
 where $\varepsilon$ is the diffusivity coefficient and $u$ is the quantity of drug injected in the domain across $\Gamma_{\text{drug}}$, which is the control variable of the problem. 
The objective functional $J$ involves the sum of two contributions, of which the first one $J_u$ related to the energy of the control, 
and the second one $J_{y_{T_f}}$ to the mismatch between the drug concentration and the desired distribution $y_{T_f}$ at the final time of the simulation $t = T_f$: 
\begin{eqnarray}\label{eq:objective-cardio}
  J(y, u)=J_u  + J_{y_{T_f}}=\frac{\beta_1}{2} \int_0^{T_f}\int_{\Gamma_{\text{drug}}} u^2 d\Omega \  dt+ 
  \frac{\beta_3}{2} \int_\Omega (y(t=T_f)-y_{T_f})^2 \ d\Omega,
\end{eqnarray}
where $\beta_1 > 0$ and $\beta_3 > 0$ are proper weights. In this case, the OCP reads: Find $u$ such that the functional \eqref{eq:objective-cardio} is minimized under the constraint that the state variable $y$ satisfies the system \eqref{eq:general_transport}. 
We define the Lagrangian functional $\mathcal{L}$ as 
\begin{align}\label{eq:lagrangian-cardio}
  \mathcal{L}(&y, u, \lambda, m_{0}, m_{\text{c}},m_{\text{inlet}},m_{\text{N}}) \nonumber\\
  &=J(y, u)-\int_0^{T_f} \int_\Omega \lambda\ (\partial_t y+\mathbf{V} \cdot \nabla y- \varepsilon\ \nabla^2 y)\ d\Omega\ dt\nonumber\\
  &+\int_\Omega m_{0}\ (y(t=0)-y_0)\ d\Omega
  + \int_0^{T_f} \int_{\Gamma_{\text{c}}} m_{\text{c}}\ (y-y_{\text{c}}) \ d\Omega\ dt\nonumber\\
  &+ \int_0^{T_f} \int_{\Gamma_{\text{inlet}}} m_{\text{inlet}}\ (y-y_{\text{in}}) \ d\Omega \ dt+\int_0^{T_f} \int_{{\Gamma_{\text{outlet}}}\cup {\Gamma_{\text{wall}}}} m_{\text{N}}\ (\nabla y\cdot \textbf{n}-y_{\text{N}}) \ d\Omega \ dt.
\end{align} 

In equation \eqref{eq:lagrangian-cardio}, $\lambda$ is the adjoint variable associated to the state variable $y$. Additionally, we also have the adjoint variables,  $m_{0}$, and $m_{\text{c}}$, $m_{\text{inlet}}$, and $m_{\text{N}}$, related to the initial data and boundary conditions. The optimal solution corresponds to
\begin{eqnarray}
\nabla \mathcal{L}(y,u, \lambda, m_{0}, m_{\text{c}}, m_{\text{inlet}}, m_{\text{N}})=(\mathcal{L}_y,\mathcal{L}_{u},\mathcal{L}_{\lambda},\mathcal{L}_{  m_{0}},\mathcal{L}_{  m_{\textbf{c}}},\mathcal{L}_{  m_{\text{inlet}}},\mathcal{L}_{  m_{\text{N}}})=\mathbf{0}.
\end{eqnarray}
By taking $\mathcal{L}_{\lambda} = 0$ and $\mathcal{L}_{m_{0}}=\mathcal{L}_{m_{\textbf{c}}}=\mathcal{L}_{m_{\text{inlet}}}=\mathcal{L}_{m_{\text{N}}} = 0$, we get the state equation with initial data and boundary conditions \eqref{eq:general_transport}. Additionally, by taking $\mathcal{L}_{y} = 0$, we get the adjoint system 
\begin{equation}
    \begin{cases}
   \partial_t \lambda+\mathbf{V}\cdot \nabla\lambda+\varepsilon \ \nabla^2 \lambda=0  & \quad \text{in} \quad \Omega \times [T_f,0), \\
   \lambda=0  & \quad \text{on} \quad \Gamma_{\text{drug}}\cup \Gamma_{\text{c}}\cup \Gamma_{\text{inlet}}\times [T_f,0),\\
   \lambda \mathbf{V}\cdot\mathbf{n}+ \varepsilon\ \nabla\lambda \cdot \mathbf{n}=0  & \quad \text{on} \quad \Gamma_{\text{outlet}} \cup \Gamma_{\text{wall}} \times [T_f,0),\\
 \lambda=\beta_3\ (y(t=T_f)-y_{T_f})  & \quad \text{in} \quad \Omega \times \{T_f\}, \\
    \end{cases} \label{eq:adjoint-cardio}
\end{equation}
and finally, by enforcing $\mathcal{L}_{u} = 0$, we obtain the optimality condition 
\begin{equation}\label{eq:optima-cardio}
     \beta_1 u - \varepsilon\ \nabla\lambda\cdot\textbf{n}-\lambda \mathbf{V}\cdot \textbf{n}=0  \quad \text{on} \quad \Gamma_{\text{drug}} \times (0, T_f]. \\
\end{equation}

The OCP is given by eqs. \eqref{eq:general_transport}, \eqref{eq:adjoint-cardio}, and \eqref{eq:optima-cardio}.

\section{Numerical discretization}\label{sec:numerical-method}
In this section we briefly discuss the time and space discretization adopted for the OCPs described in the previous section. 

We start with the time discretization. For what concerns the state problem, we choose a time step $\Delta t \in \mathbb{R}^+$ to subdivide the time interval $(0,T_f]$ into $t^n =  n \Delta t$, with $n = 0, ..., N_{tf}$ and $T_f =  N_{tf} \Delta t$. On the other hand, for the adjoint problem, we switch the sign of the time derivative, i.e., we perform a reflection in time in order to consider the same collocation points as in the state problem. 
The time derivatives are discretized by means of the Euler method. Both convective and diffusive terms are treated implicitly. 

For the space discretization, the computational domain $\Omega$ is divided into cells or control volumes $\Omega_i$, with $i = 1, \dots, N_{c}$, where $N_{c}$ is the total number of cells in the mesh. We adopt a second-order central finite volume scheme 
both for convective and diffusive terms. 

We opt for a segregated approach to solve our OCPs. See Algorithm \ref{alg:gradient_opt}. 

\begin{algorithm}[H]
\caption{}
\begin{algorithmic}[1]
\State Provide an \textbf{initial guess} for the control variable;
\Repeat
    \State Solve the \textbf{state problem};
    \State Solve the \textbf{adjoint problem};
    \State \textbf{Update} the control variable by solving the optimality condition 
    using the steepest gradient method \cite{kelley1999iterative};
\Until{$|\mathcal{L}_u| < tol$}.
\end{algorithmic}
\label{alg:gradient_opt}
\end{algorithm}
In all our numerical experiments, the initial guess for the control variable is set to zero while the tolerance is $tol=10^{-6}$.

All of the above choices are so that we can design an efficient splitting method for the complex coupled problem at hand. 

The above computational pipeline is implemented in the open-source finite volume C++ library OpenFOAM\textsuperscript{\textregistered} \cite{openFoam}. 

\section{Computational results}
\label{sec:computational_results}
This section presents the numerical results of the OCPs described in Sec. \ref{ocp} obtained by using the algorithm reported in Sec. \ref{sec:numerical-method}. 
The academic benchmark is discussed in Sec. \ref{sec:analytical} while the two biomedical test cases are presented in Secs. \ref{sec:drugdelivery} and \ref{sec:drugcoronary}. 
\subsection{Academic benchmark}\label{sec:analytical}
By following \cite{fu2009} let us $\Omega$ be the unit square [$0, 1]^2$. We set \( \mathbf{V} = (1, 0) \). In order to validate our numerical approach, we use a manufactured solution:  the functions $f$ and $y_d$ are chosen such that the exact expressions for the state, adjoint, and control variables are given by: 
\begin{eqnarray}
 && y_\text{exact}= e^{-t} \sin(2\pi x_1) \sin(2\pi x_2),\nonumber\\
 &&\lambda_\text{exact}= e^{-t} (T_f - t) \sin(2\pi x_1) \sin(2\pi x_2),\nonumber\\
 && u_\text{exact} = -\frac{ e^{-t} (T_f - t)}{\beta_1} \sin(2\pi x_1) \sin(2\pi x_2). \label{manufactured-state-costate}
\end{eqnarray}
For sake of completeness and in order to make easier the reproducibility of the results, we also report the expressions of $f$ and $y_d$: 
\begin{align}
 f = e^{-t} \sin(2\pi x_1) \sin(2\pi x_2) \big(8\pi^2 \epsilon -1+\frac{1}{\beta_1} (T_f-t)\big)+2\pi \ e^{-t}\ \cos(2\pi x_1) \sin(2\pi x_2), \\
 y_d =\frac{e^{-t} (T_f - t) \sin(2\pi x_1) \sin(2\pi x_2)}{\beta_2} \big((t-1-T_f)-8\pi^2  \epsilon (T_f-t)\big)
 +  \notag \\ \frac{2\pi e^{-t} (T_f-t) \cos(2\pi x_1) \sin(2\pi x_2)}{\beta_2} +e^{-t} \sin(2\pi x_1) \sin(2\pi x_2).
\end{align}

We set $T_f = 1$ and $\beta_1 = \beta_2 = 1$. We consider mesh sizes $h^{(i)} = 1/2^{i+1}$ for $i = 1, 2, 3, 4$. The time step $\Delta t$ is chosen such that ${h^{(i)}}^2 \Delta t = 1$.  
We compute the rate of convergence as: 
\begin{eqnarray}\label{analytical-rate}
    \text{rate} = \frac{\log\left(\frac{E^{(i)}_\varphi}{E^{(i+1)}_\varphi}\right)}{\log\left(\frac{h^{(i)}}{h^{(i+1)}}\right)}, \quad \text{for} \ i = 1, 2, 3,
\end{eqnarray}
where $E^{(i)}_\varphi$ is the $L^2$-norm relative error computed on mesh $h^{(i)}$ between the numerical solution and the analytical one (see eq. \eqref{manufactured-state-costate}): 
\begin{eqnarray}\label{analytical-error}
    E^{(i)}_\varphi = \dfrac{||\varphi^{(i)} - \varphi_{\text{ex}}||_{\Omega}}{||\varphi_\text{ex}||_{\Omega}}, \quad \text{where} \ \varphi = \{y, \lambda \}.
\end{eqnarray}

Tables~\ref{table:ms-time-ind-error-ep-01}--\ref{table:ms-time-ind-error-ep-2} report the values of $E^{(i)}_\varphi$ as well as the convergence rates for the state and adjoint variables, $y$ and $\lambda$, at the final time $T_f=1$, for three different values of the diffusion coefficient $\varepsilon = 1, 10^{-1}, 10^{-2}$, respectively. 
 We also report the number of iterations $N$ required to reach convergence.   
It should be noted that the convergence rate remains consistently at or near 2.00 for both variables. 
This is expected, as a second-order approximation is used for the discretization of spatial derivatives (see Sec. \ref{sec:numerical-method}). 
Concerning the number of iterations $N$, it is almost constant with respect to the mesh size at a given value of $\varepsilon$. However, it slightly increases with decreasing $\varepsilon$. 

For further comparison, Figs. \ref{fig:manufactured_y_re_error1} and \ref{fig:manufactured_lambda_error1} display the numerical solution, the analytical one and the related pointwise error at the final time of the simulation $T_f = 1$ on the mesh $h = 1/32$. 
\begin{table}[h]
    \centering
    \resizebox{0.6\textwidth}{!}{
    \begin{tabular}{|c|c|c|c|c|c|}
       \hline
       \multirow{2}{*}{$h$} & \multirow{2}{*}{$N$} & \multicolumn{2}{|c|}{$E_y$} & \multicolumn{2}{|c|}{$E_\lambda$} \\ \cline{3-6} 
        & &  error & rate & error & rate   \\ \hline
        $1/4$ & 3  & $2.3798 \cdot 10^{-1}$ & / & $1.5702 \cdot 10^{-1}$ &  /  \\ \hline
        $1/8$ & 3 & $5.3124 \cdot 10^{-2}$ & 2.16 & $3.5380 \cdot 10^{-2}$ & 2.15   \\ \hline
        $1/16$ & 3 & $1.2321 \cdot 10^{-2}$ & 2.10 & $8.0509 \cdot 10^{-3}$ & 2.13   \\ \hline
        $1/32$ & 3 & $2.8101 \cdot 10^{-3}$ & 2.13 & $2.0358 \cdot 10^{-3}$ & 1.98   \\ \hline
    \end{tabular}}
   \caption{Academic benchmark: number of iterations $N$, $L^2$-norm relative error between analytical and numerical solution (eq. \eqref{analytical-error}) and rates of convergence  (eq. \eqref{analytical-rate}) for $\epsilon = 1$ computed on different mesh sizes.}
    \label{table:ms-time-ind-error-ep-01}
\end{table}


\begin{table}[h]
    \centering
    \resizebox{0.6\textwidth}{!}{
    \begin{tabular}{|c|c|c|c|c|c|}
       \hline
       \multirow{2}{*}{$h$} & \multirow{2}{*}{$N$} & \multicolumn{2}{|c|}{$E_y$} & \multicolumn{2}{|c|}{$E_\lambda$} \\ \cline{3-6} 
        & & error & rate & error & rate   \\ \hline
        $1/4$ & 5 & $3.8452 \cdot 10^{-1}$ & / & $1.9788 \cdot 10^{-1}$ &  /  \\ \hline
        $1/8$ & 5 & $8.1619 \cdot 10^{-2}$ & 2.23 & $4.5045 \cdot 10^{-2}$ & 2.13   \\ \hline
        $1/16$ & 5 & $1.9373 \cdot 10^{-2}$ & 2.07 & $1.0567 \cdot 10^{-2}$ & 2.09   \\ \hline
        $1/32$ & 5 & $4.7043 \cdot 10^{-3}$ & 2.04 & $2.7235 \cdot 10^{-3}$ & 1.95   \\ \hline
    \end{tabular}}
   \caption{Academic benchmark: number of iterations $N$, $L^2$ relative error between analytical and numerical solution (eq. \eqref{analytical-error}) and rates of convergence (eq. \eqref{analytical-rate}) for $\epsilon = 10^{-1}$ computed on different mesh sizes.}
    \label{table:ms-time-ind-error-ep-1}
\end{table}


\begin{table}[h]
    \centering
    \resizebox{0.6\textwidth}{!}{
    \begin{tabular}{|c|c|c|c|c|c|}
       \hline
       \multirow{2}{*}{$h$} & \multirow{2}{*}{$N$} & \multicolumn{2}{|c|}{$E_y$} & \multicolumn{2}{|c|}{$E_\lambda$} \\ \cline{3-6} 
        & & error & rate & error & rate   \\ \hline
        $1/4$ & 10 & 1.3322 &  /& $2.0563 \cdot 10^{-1}$ & /   \\ \hline
        $1/8$ & 10 & $1.9883 \cdot 10^{-1}$ & 2.74 & $5.1639 \cdot 10^{-2}$ & 1.99   \\ \hline
        $1/16$ & 10 & $3.9698 \cdot 10^{-2}$ & 2.32 & $1.1536 \cdot 10^{-2}$ & 2.16   \\ \hline
        $1/32$ & 9 & $7.2226 \cdot 10^{-3}$ & 2.45 & $2.5999 \cdot 10^{-3}$ & 2.14   
        \\ \hline
    \end{tabular}}
   \caption{Academic benchmark: number of iterations $N$, $L^2$ relative error between analytical and numerical solution (eq. \eqref{analytical-error}) and rates of convergence (eq. \eqref{analytical-rate}) for $\epsilon = 10^{-2}$ computed on different mesh sizes.}
    \label{table:ms-time-ind-error-ep-2}
\end{table}

\begin{figure}[htb!]
\centering
\begin{tabular}{ccccc}
     & $\epsilon=1$ & $\epsilon=10^{-1}$ & $\epsilon=10^{-2}$  & \\
\raisebox{1.5cm}{$y$} & \includegraphics[scale = 0.1]{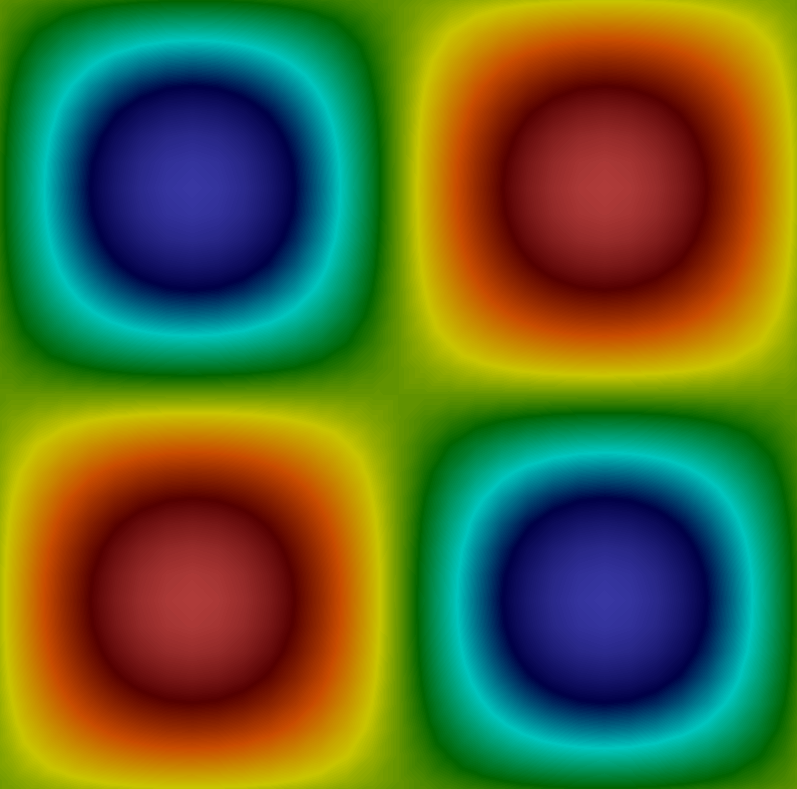} & \includegraphics[scale = 0.1]{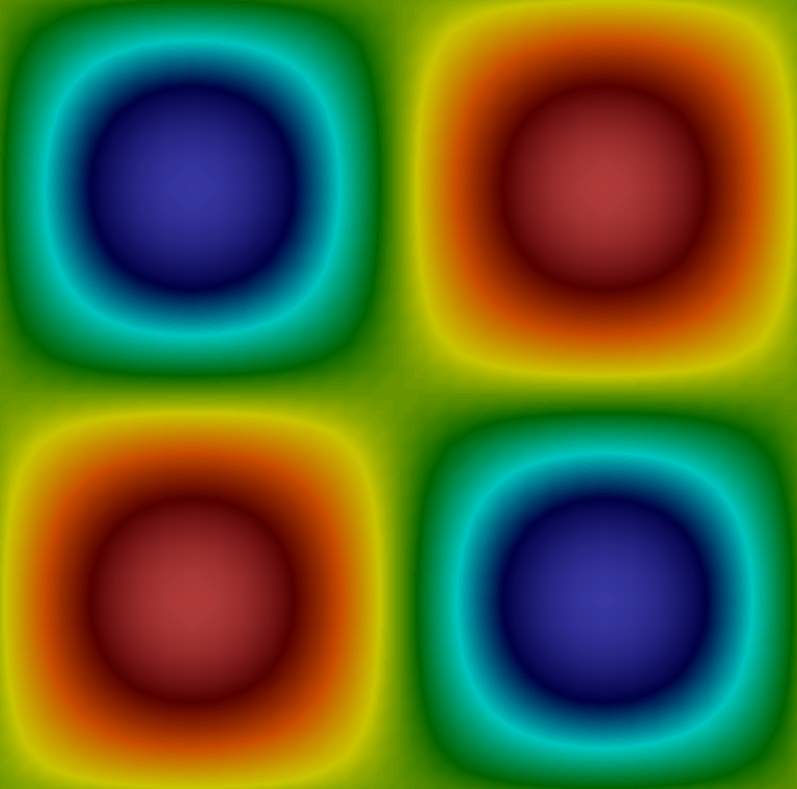} & \includegraphics[scale = 0.1]{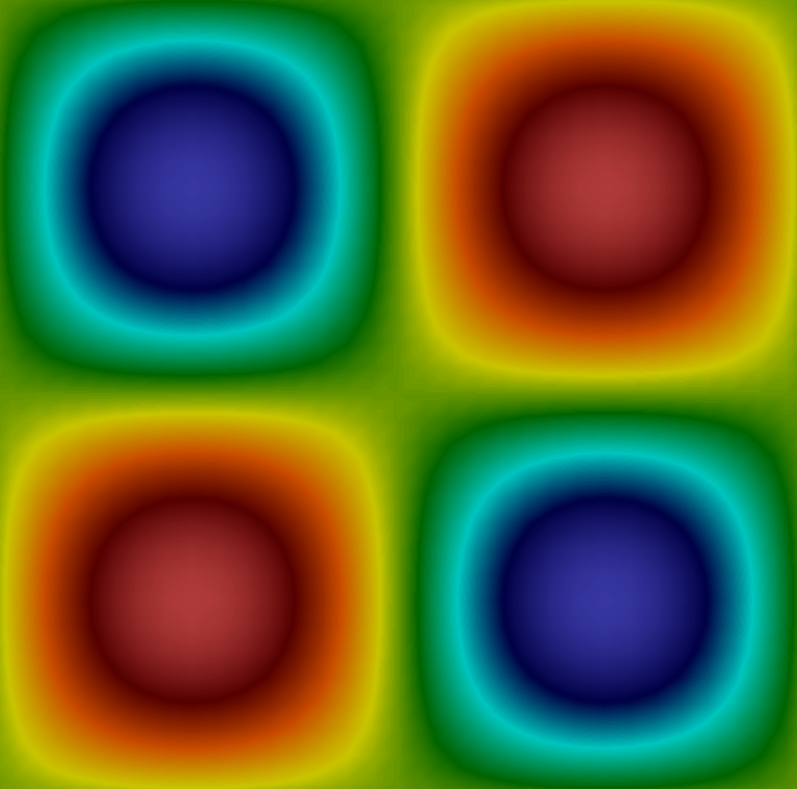} & 
\hspace{-0.4cm} 
\includegraphics[width=1.5cm, height=3cm]{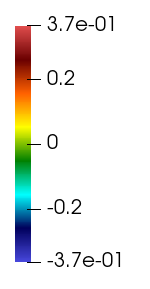} \\
\raisebox{1.5cm}{$y_\text{exact}$}  & \includegraphics[scale = 0.1]{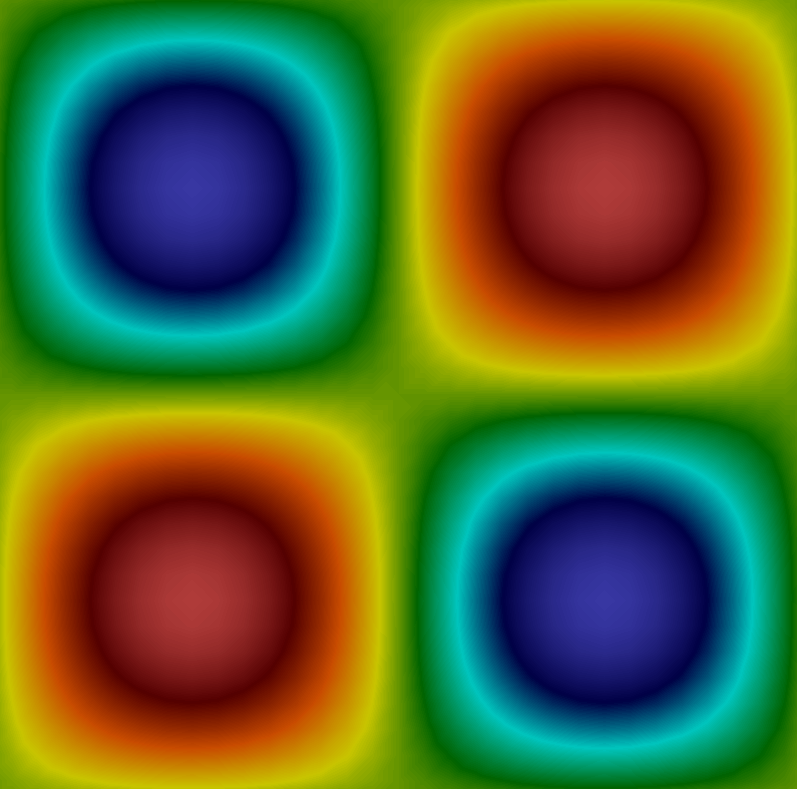} & \includegraphics[scale = 0.1]{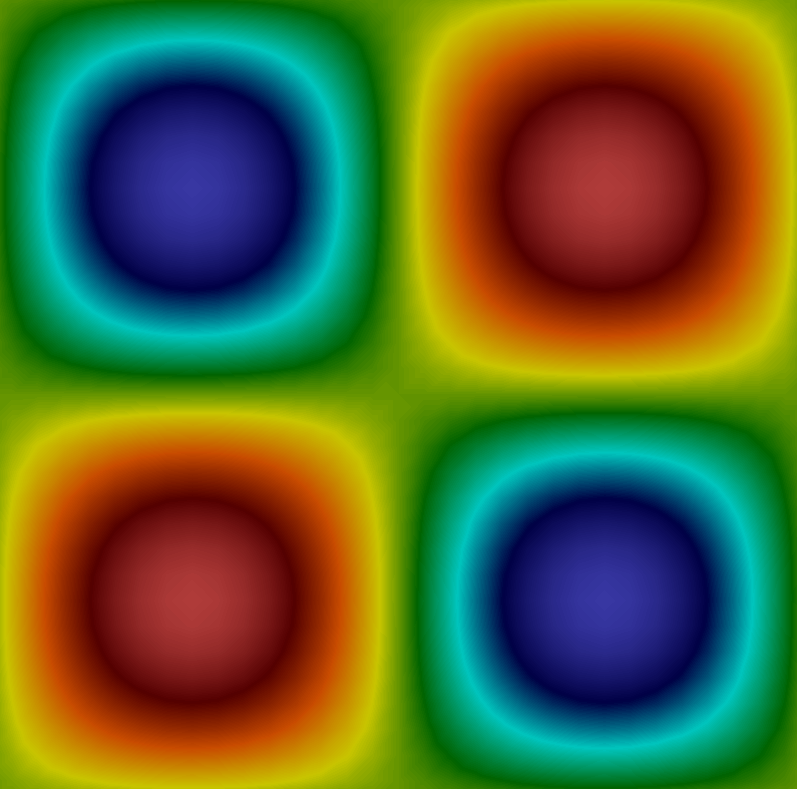} & \includegraphics[scale = 0.1]{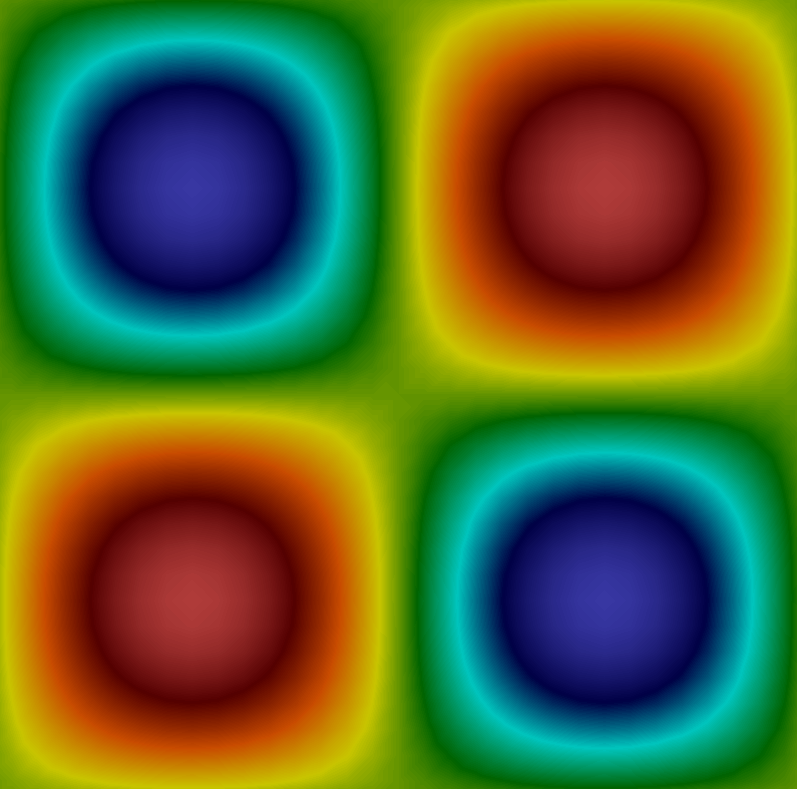} 
& \hspace{-0.4cm} \includegraphics[ width=1.5cm, height=3cm]{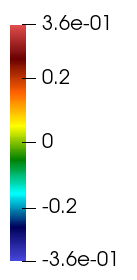}\\
\raisebox{1.5cm}{$y-y_\text{exact}$}  & \includegraphics[scale = 0.1]{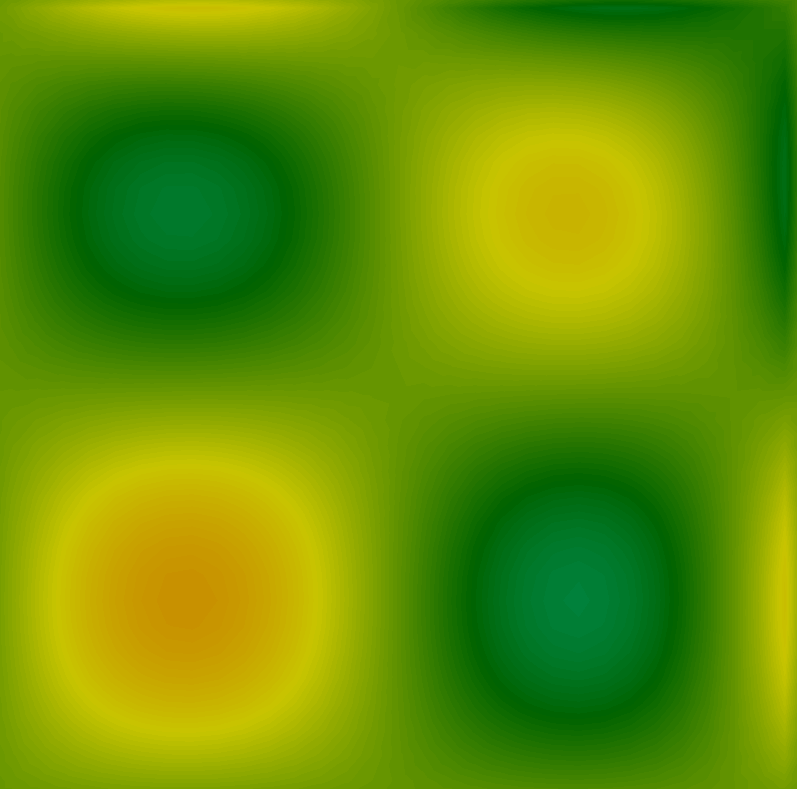} & \includegraphics[scale = 0.1]{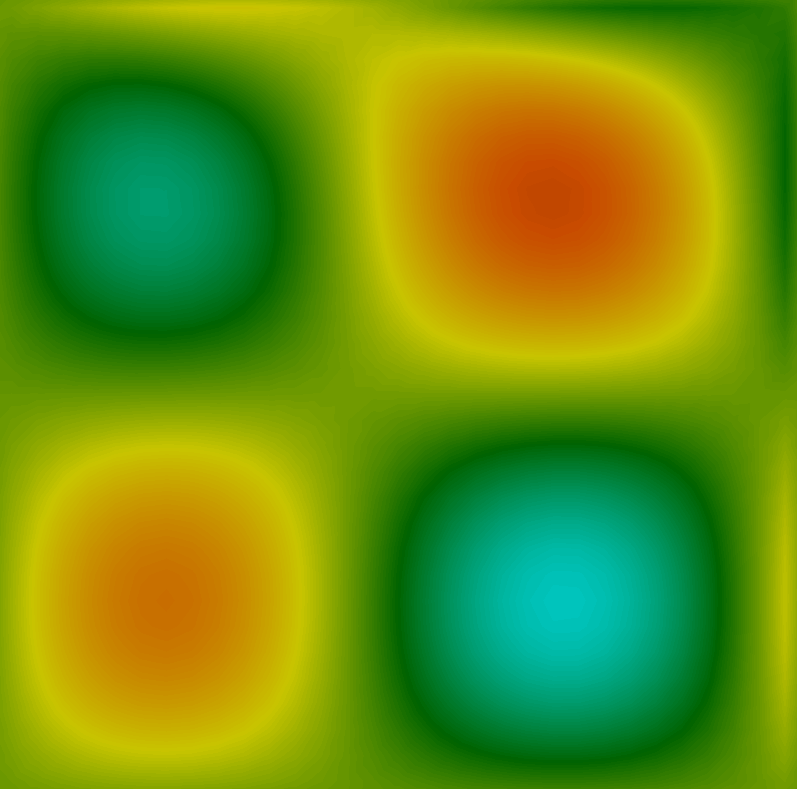} & \includegraphics[scale = 0.1]{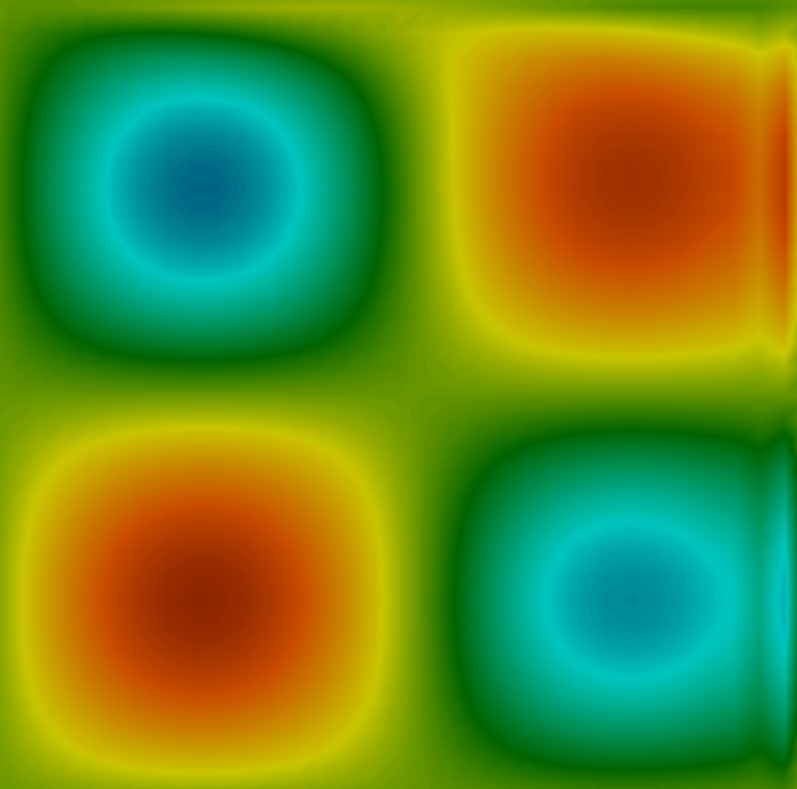} 
&  \hspace{-0.4cm} \includegraphics[width=1.5cm, height=3cm]{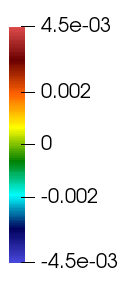}
\end{tabular}
\caption{Academic benchmark: numerical solution $y$ (top row), analytical solution $y_{\text{exact}}$ (second row) and pointwise error $y-y_\text{exact}$ (bottom row) for $\epsilon = 1$ (first column), $\epsilon = 10^{-1}$ (second column), and $\epsilon = 10^{-2}$ (third column) at the final time $T_f=1$. The results were obtained with a mesh of size $h = 1/32$. } 
\label{fig:manufactured_y_re_error1}
\end{figure}

\begin{figure}[htb!]
\centering
\begin{tabular}{ccccc}
      & $\epsilon=1$ & $\epsilon=10^{-1}$ & $\epsilon=10^{-2}$  & \\
\raisebox{1.5cm}{$\lambda$}  & \includegraphics[scale = 0.1]{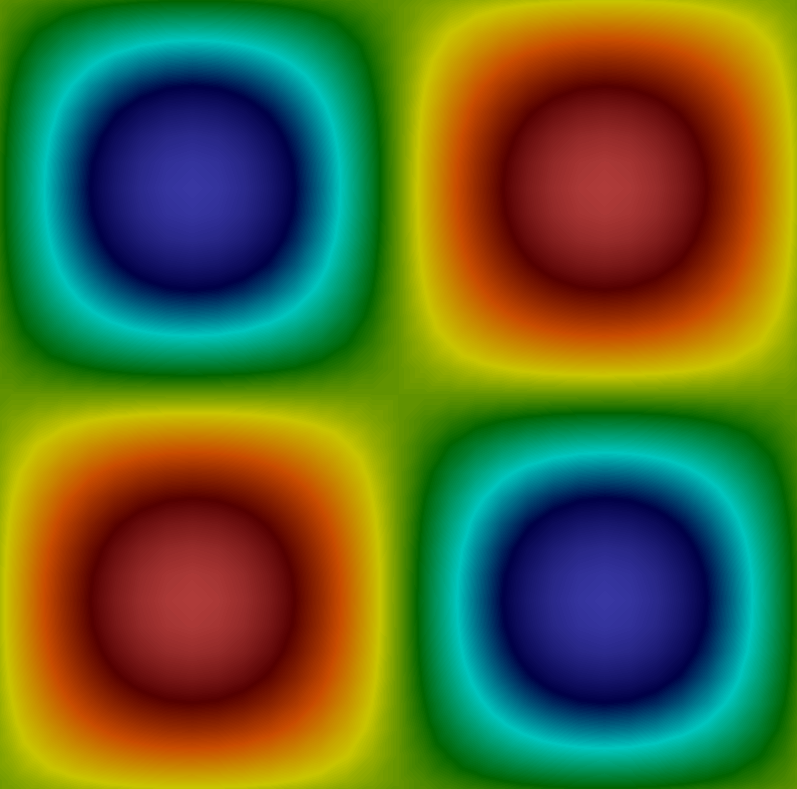} & \includegraphics[scale = 0.1]{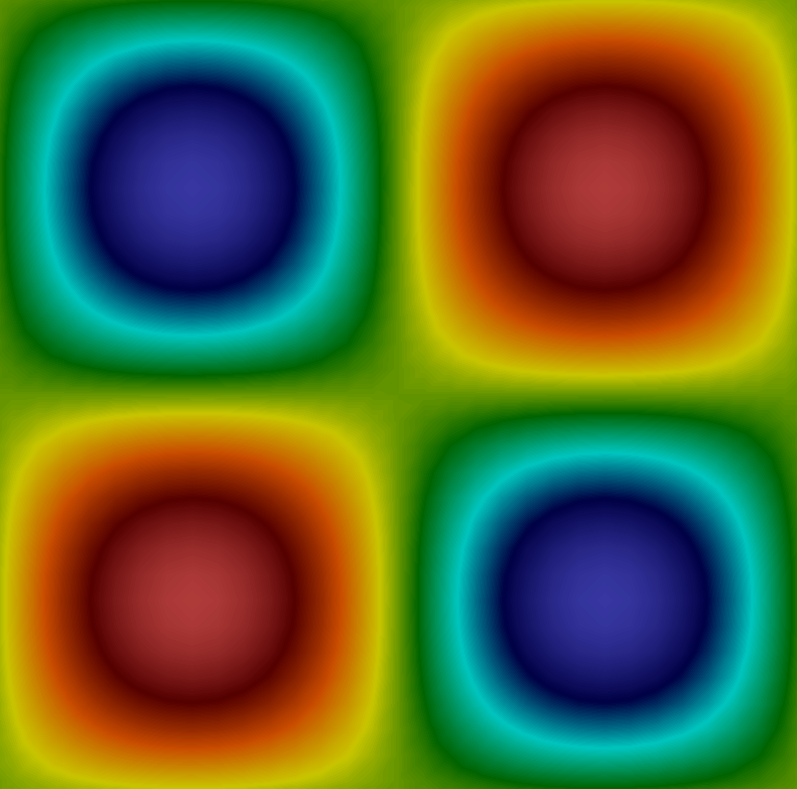} & \includegraphics[scale = 0.1]{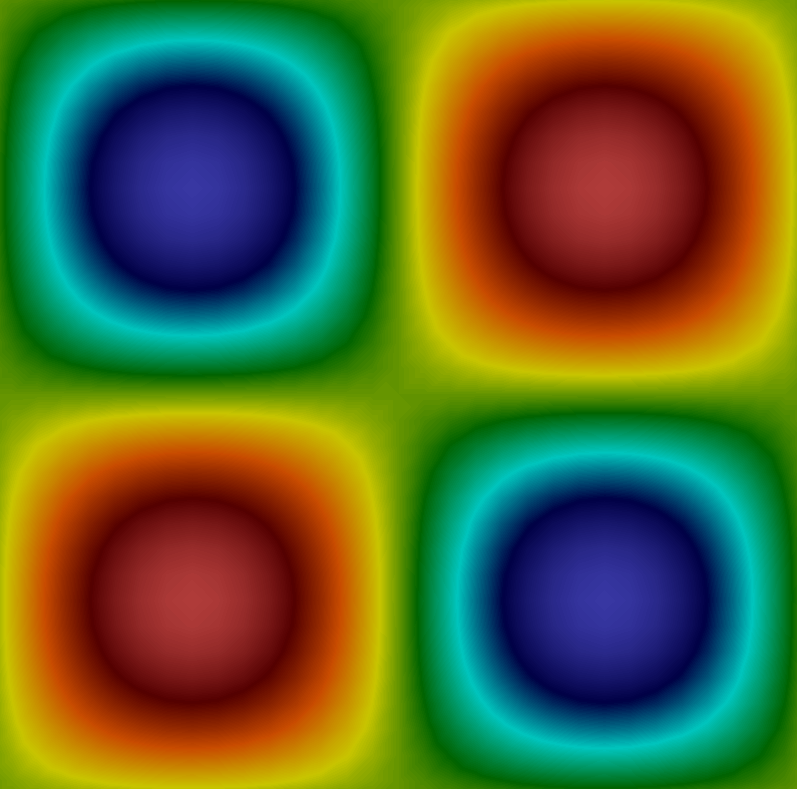} 
& \hspace{-0.4cm} \includegraphics[width=1.5cm, height=3cm]{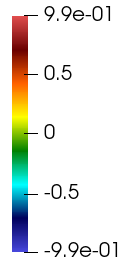} \\
\raisebox{1.5cm}{$\lambda_\text{exact}$}  & \includegraphics[scale = 0.1]{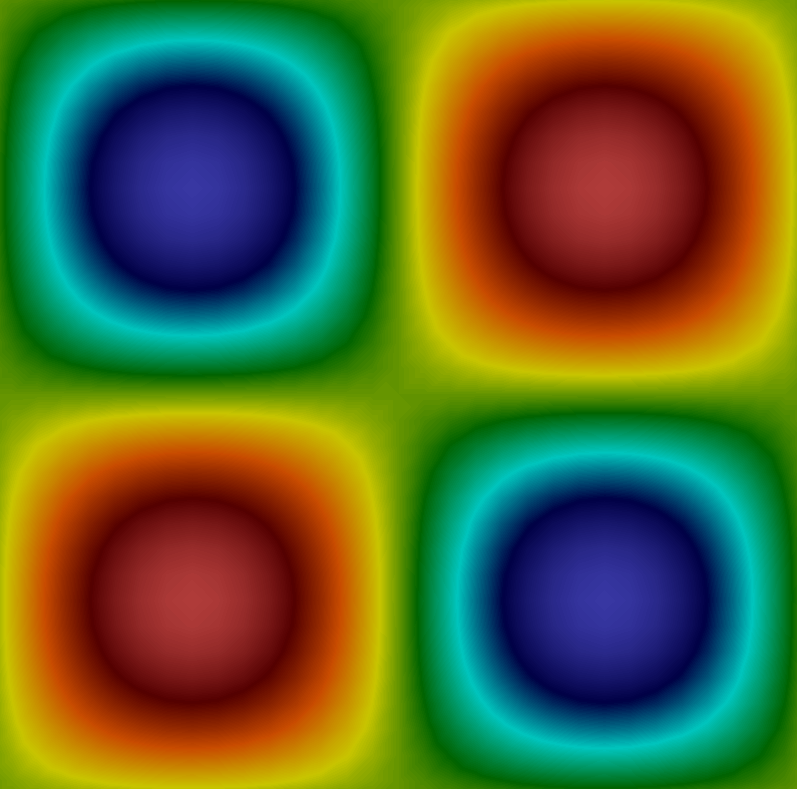} & \includegraphics[scale = 0.1]{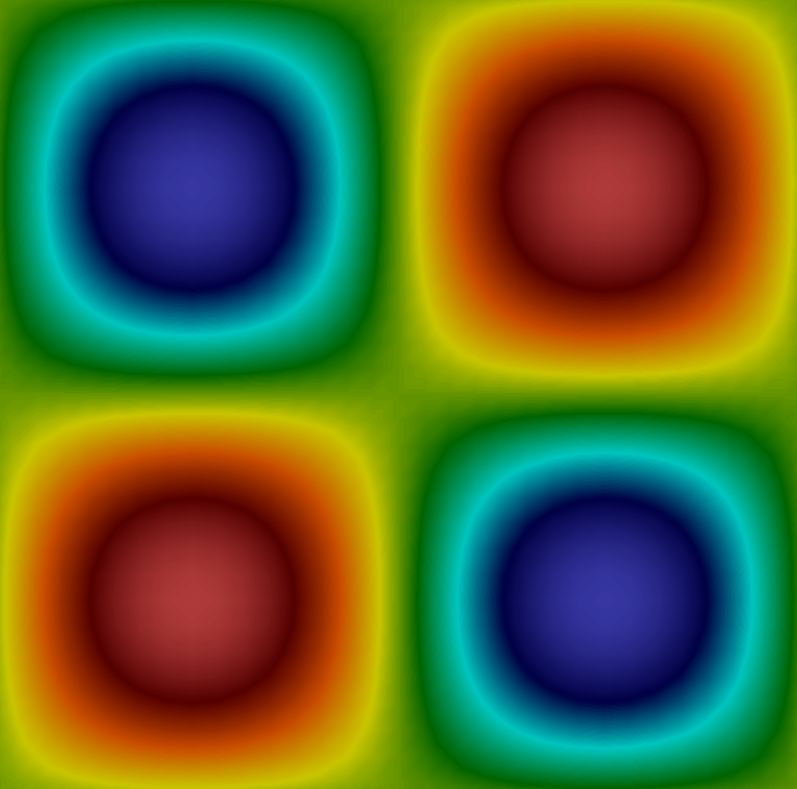} & \includegraphics[scale = 0.1]{img_manufactured_lambda_anal_epsilon1e-2_32.png} 
& \hspace{-0.4cm} \includegraphics[width=1.5cm, height=3cm]{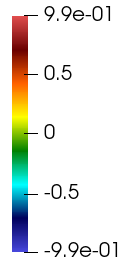}\\
\raisebox{1.5cm}{$\lambda-\lambda_\text{exact}$}  & \includegraphics[scale = 0.1]{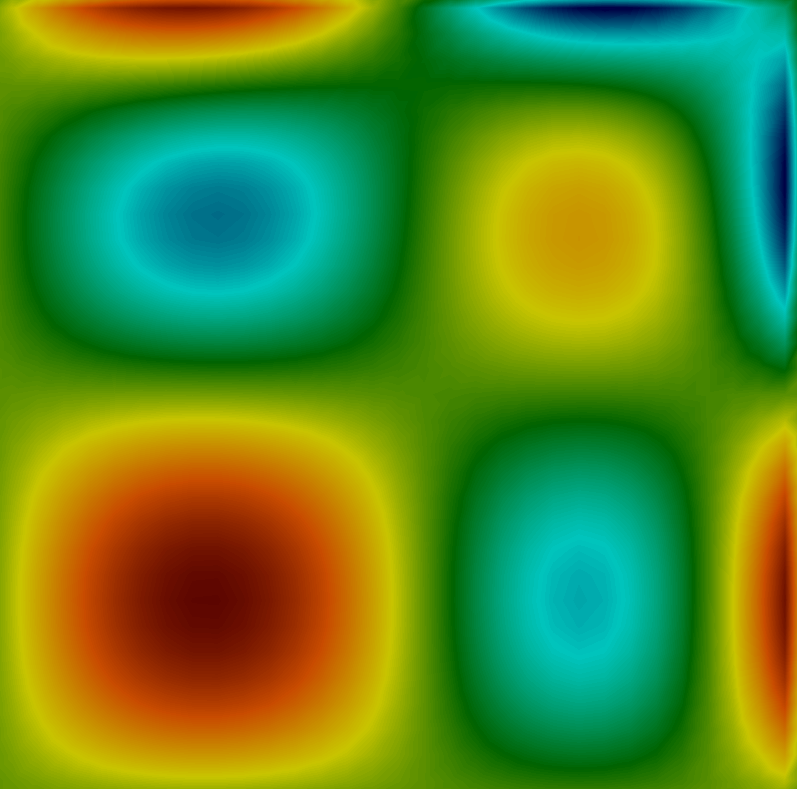} & \includegraphics[scale = 0.1]{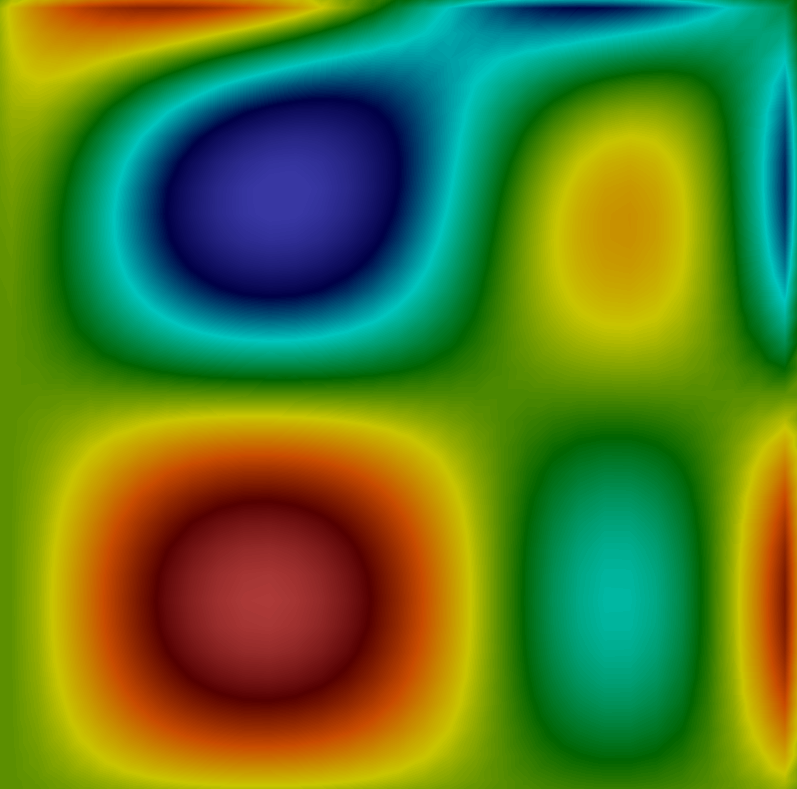} & \includegraphics[scale = 0.1]{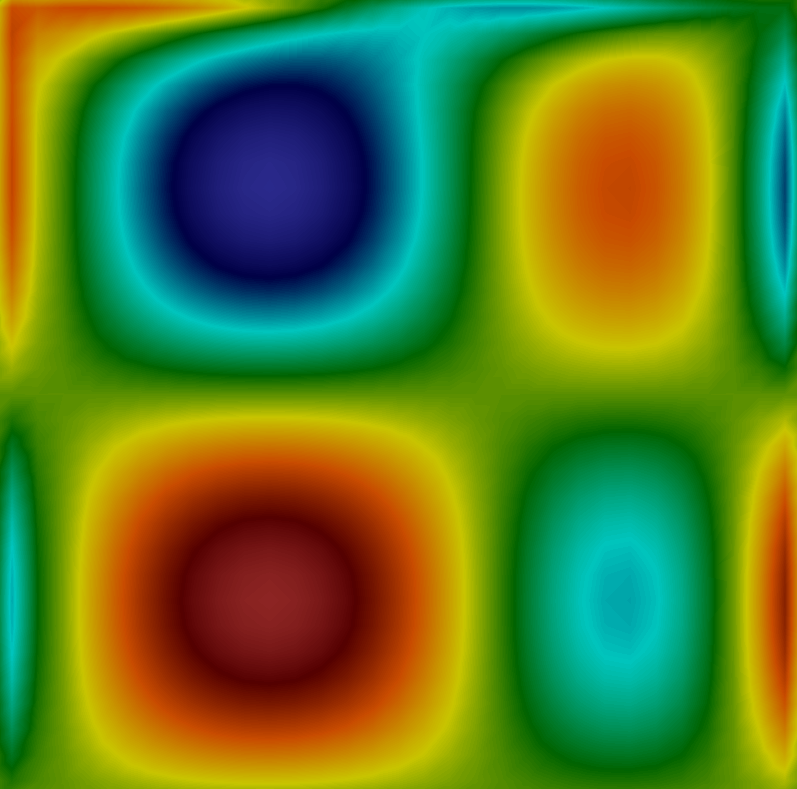} 
&  \hspace{-0.4cm} \includegraphics[width=1.5cm, height=3cm]{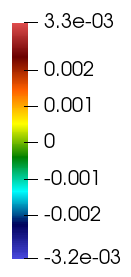}
\end{tabular}
\caption{Academic benchmark: numerical solution $\lambda$ (top row), analytical solution $\lambda_{\text{exact}}$ (second row) and pointwise error $\lambda-\lambda_\text{exact}$ (bottom row) for $\epsilon = 1$ (first column), $\epsilon = 10^{-1}$ (second column), and $\epsilon = 10^{-2}$ (third column) at the final time $T_f=1$. The results were obtained with mesh size $h = 1/32$.}
\label{fig:manufactured_lambda_error1}
\end{figure}

\subsection{Application 1: Drug delivery in a light-driven system} \label{sec:drugdelivery}

\subsubsection{Distributed case}

By following \cite{ferreira2023prepub} 
we set \(\Omega = [0,1]\). In addition, we set $D_d=4 \times 10^{-4}$ and $\gamma = 4 \times 10^{-3}$.  
For initial data and boundary conditions we have: $c_{f,0} = 0$, $c_{b,0} = 1$ for $0 \leq x \leq 0.25$ and  $c_{b,0} = 0$ for $0.25 < x \leq  1$, $c_{f,r} = 0$ and $c_{fN} = 0$. The mesh size is $h = 256$. We run all the simulations until $T_f = 10$ with $\Delta t = 0.1$. 

In order to validate our numerical approach, since we do not have any clinically relevant free drug profile to be assimilated in our model, we generate fictitious data by solving the state problem \eqref{light-unsteady-enhancer} with 
a certain value of $u$. In particular, we employ the configurations tested in \cite{ferreira2023prepub}, i.e., $u = I^0 e^{-\omega x}$ with $\omega= 4$ and two values of $I^0 \in \{5,15\}$. Then 
we consider as target data $c_{f,T_f}$ the solution $c_f (t = T_f)$ at the final time of the state problem simulation and we solve the OCP problem \eqref{light-unsteady-enhancer}, \eqref{adjointdebc-unsteady-light-distributed-lambda1-lambda2} and \eqref{eq:optimality}. We set $\beta_3 = 1$ 
while we perform a sensitivity analysis with respect to $\beta_1 \in \{10^{-3}, 10^{-4}, 10^{-5}, 10^{-6}\}$. 



The results are shown in Figures ~\ref{fig:plots_light_1D_I05_distributed}--\ref{fig:plots_light_1D_I15_distributed}. 
As expected, as $\beta_1$ decreases,  both $J$ and $J_{c_{f,T_f}}$ assume smaller and smaller values, although convergence is reached after a greater number of iterations. 
Consequently, the OCP solution provides a spatial distribution of the control variable $u$ which is closer to $I^0 e^{-\omega x}$ (i.e., the source term which we have employed in the state problem to generate the target data) for lower values of $\beta_1$. In particular, we note that for $\beta_1 = 10^{-6}$, except for a slight mismatch next to $x = 0 $ and $x = 0.25$, $u$ and $I^0 e^{-\omega x}$ are in very good agreement. We highlight that the comparison is shown only until $x = 0.25$ because, as the control variable $u$ is proportional to $c_b$ (see eq. \eqref{eq:optimality}) 
and $c_b$ is null beyond $x = 0.25$, the control variable is also null beyond $x = 0.25$. 
Finally, the goodness of the solution is also demonstrated by the comparison between the free drug concentration at the final time of the OCP simulation (i.e., when the control variable is adopted as a source term) and the target value $c_{f, T_f}$. 

\begin{figure}[htb!]
    \centering
    \begin{subfigure}[b]{0.48\textwidth}
        \centering
        \includegraphics[scale=0.48]{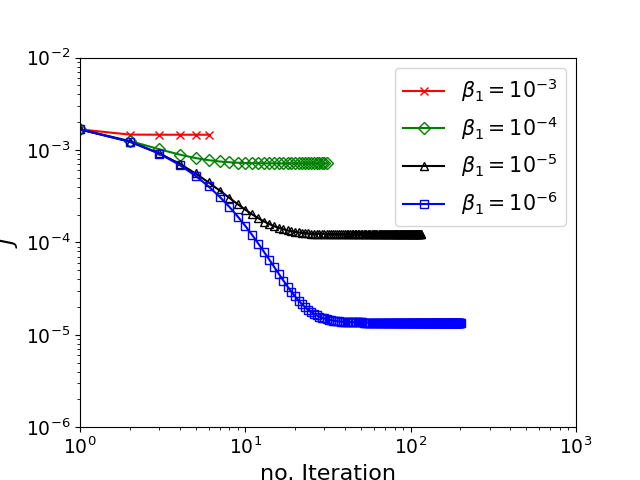}
    \end{subfigure}
    \begin{subfigure}[b]{0.48\textwidth}
        \centering
        \includegraphics[scale=0.48]{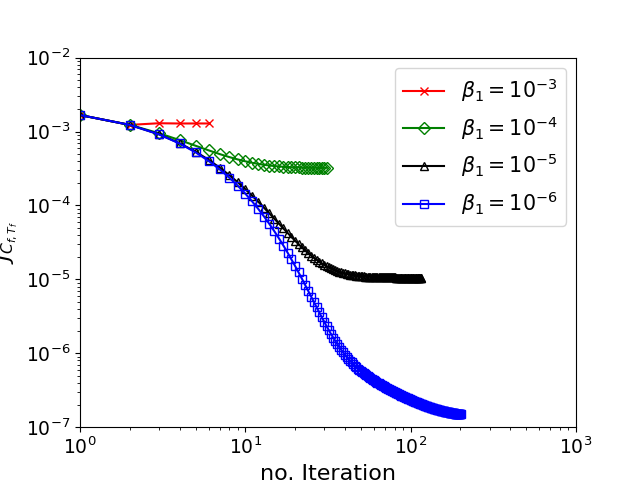}
    \end{subfigure}
    \\ 
    \begin{subfigure}[b]{0.48\textwidth}
        \centering
        \includegraphics[scale=0.48]{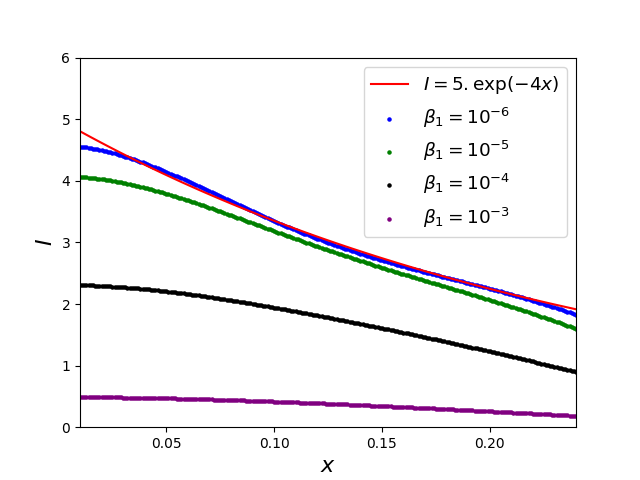}
    \end{subfigure}
    \begin{subfigure}[b]{0.48\textwidth}
        \centering
        \includegraphics[scale=0.48]{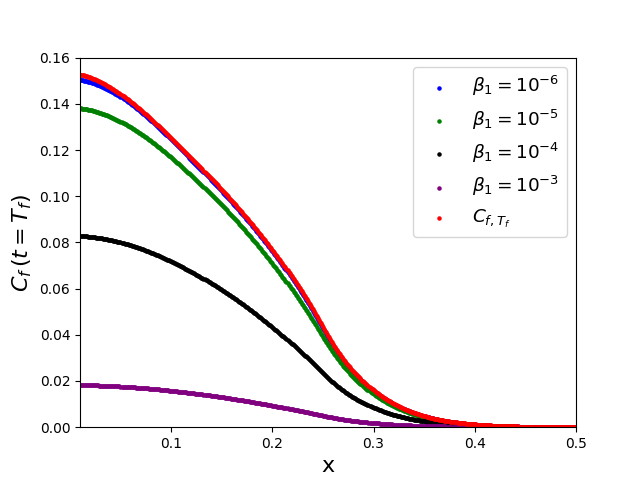}
    \end{subfigure}
    \caption{Application 1 - Distributed case: convergence history of the objective function $J$ (top left) and of the contribution $J_{c_{f,T_f}}$ (top right), spatial distribution of the control variable (bottom left) and of the free drug concentration at the final time of the simulation (bottom right) for different values of $\beta_1$ and for $I^0=5$.}
    \label{fig:plots_light_1D_I05_distributed}
\end{figure}


\begin{figure}[htb!]
    \centering
    \begin{subfigure}[b]{0.48\textwidth}
        \centering
        \includegraphics[scale=0.48]{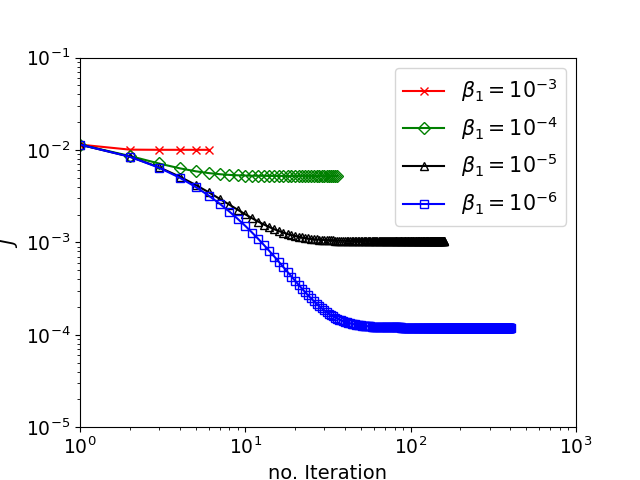}
    \end{subfigure}
    \begin{subfigure}[b]{0.48\textwidth}
        \centering
        \includegraphics[scale=0.48]{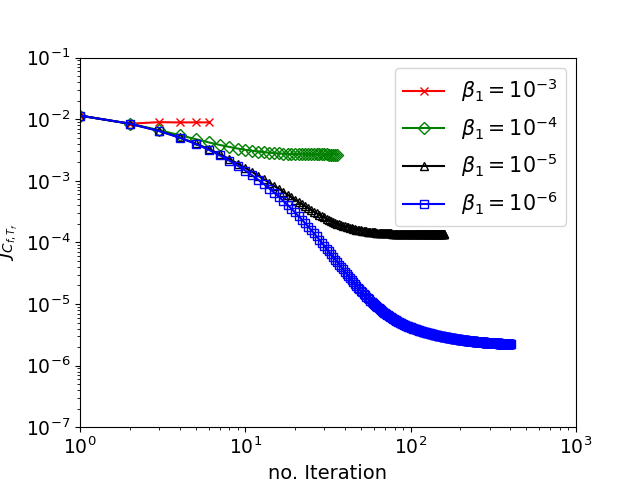}
    \end{subfigure}
    \\ 
    \begin{subfigure}[b]{0.48\textwidth}
        \centering
        \includegraphics[scale=0.48]{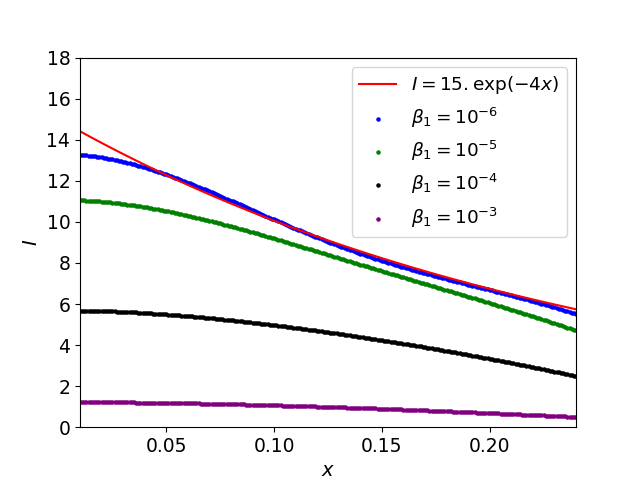}
    \end{subfigure}
    \begin{subfigure}[b]{0.48\textwidth}
        \centering
        \includegraphics[scale=0.48]{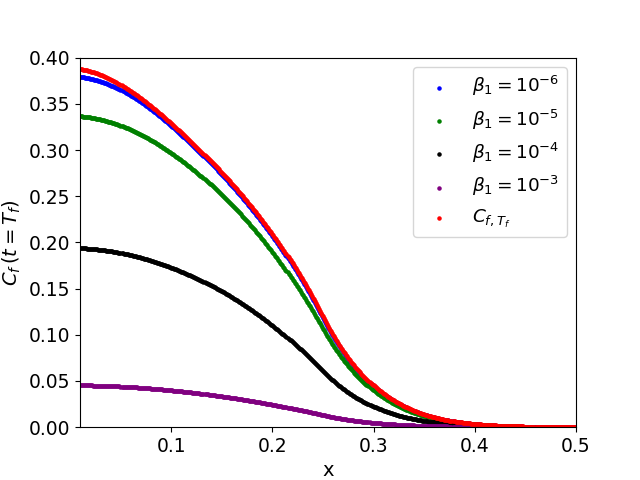}
    \end{subfigure}
    \caption{Application 1 - Distributed case: convergence history of the objective function $J$ (top left) and of the contribution $J_{c_{f,T_f}}$ (top right), spatial distribution of the control variable (bottom left) and of the free drug concentration at the final time of the simulation (bottom right) for different values of $\beta_1$ and for $I^0=15$.}
    \label{fig:plots_light_1D_I15_distributed}
\end{figure}

\subsubsection{Concentrated case}  


As reported at the beginning of Sec. \ref{sec:2_2} 
we consider both one and two-dimensional configurations. 

In the one dimensional case, 
we set \(\Omega = [0,1]\) \cite{ferreira2023} while in the two dimensional case \(\Omega = [0,1]^2\). In addition, we set $D_I=4\times 10^{-3}$, $D_d=4 \times 10^{-4}$, $\mu_a=4\times 10^{-3}$%
 and $\gamma = 15 \times 10^{-3}$. 
For initial data and boundary conditions we have: $I_N = 0$, $I_0 = 0$, $c_{f,0} = 0$, $c_{b,0} = 1$ for $0 \leq x \leq 0.25$ and  $c_{b,0} = 0$ for $0.25 < x \leq  1$, $c_{f,r} = 0$ and $c_{fN} = 0$. The mesh consists of 256 cells for the one dimensional configuration whilst of 4096 cells for the two dimensional one. We run all the simulations until $T_f = 10$ with $\Delta t = 0.1$. 

The validation is performed by using the same pipeline adopted for the distributed case. We generate data by solving the state problem \eqref{light-concentrated-unsteady-enhancer} with an assigned value of $u$. For the one-dimensional case, we employ the configurations tested in \cite{ferreira2023}, i.e., $u = I^0$.  For the two-dimensional case, we employ $u = I^0 y(1-y)$. In both cases, $I^0 \in \{5,15\}$. Then 
we consider as target data $c_{f,T_f}$ the solution $c_f (t = T_f)$ at the final time of the state problem simulation and  we solve the OCP problem  \eqref{light-concentrated-unsteady-enhancer}, \eqref{adjoint-unsteady-light-lambda1-lambda2-lambda3} and \eqref{eq:light-optimality-concentrated}. Like in the distributed case, we set $\beta_3 = 1$ 
while we perform a sensitivity analysis with respect to $\beta_1 \in \{10^{-3}, 10^{-4}, 10^{-5}, 10^{-6}\}$. 



The results are shown in Figures \ref{fig:plots_light_1D_I05_concentrated}-\ref{fig:plots_light_1D_I15_concentrated} for the one dimensional case and in Figures \ref{fig:plots_light_2D_I05_concentrated}-\ref{fig:plots_light_2D_I015_concentrated_error}
for the two dimensional case. The same observations made for the distributed case still hold. As $\beta_1$ decreases,  both $J$ and 
$J_{c_{f,T_f}}$ assume smaller and smaller values and the convergence is reached after a greater number of iterations. 
In the two dimensional case we note that the OCP solution provides a spatial distribution of the control variable which is closer to $I^0 \ y \ (1-y)$ (i.e., the boundary conditions enforced on the left boundary of the domain in the state problem to generate the target data) for lower values of $\beta_1$. In particular we note that for $\beta_1 = 10^{-6}$, unless a slight mismatch at the initial part and the final part of the domain, the two fields are in very good agreement. On the other hand, in the one dimensional case,  we do not have any spatial distribution of the control variables, i.e., it degenerates to a single value. So, in order to verify the accuracy of the OCP solution, we report in Table~\ref{table:summary_Il_beta1}  
the relative error defined as follows:
\begin{align}\label{eq:err2}
E_{u} = \frac{|u - I^0|}{I^0}.
\end{align}
We observe that the error is less than 3\% when $\beta_1 = 10^{-6}$ for all the tested configurations.
In conclusion, we comment on the comparison between the free drug concentration at the final simulation time predicted by the OCP problem (i.e., with the control variable applied as the left boundary condition) and the target distribution. The match 
is excellent for both the one- and two-dimensional configurations.  




\begin{table}[h!]
\centering
\begin{tabular}{|c|c|c|c|}
\hline
$I^0$ & $\beta_1$ 
& $u$ & $E_u$ \\
\hline
5  & $10^{-3}$ 
& 0.45  & 90.94\% \\
  & $10^{-4}$ 
  & 2.36  & 52.74\% \\
  & $10^{-5}$ 
  & 4.44  & 11.26\% \\
  & $10^{-6}$ 
  & 4.92  & 1.50\% \\
\hline
15 & $10^{-3}$ 
& 1.17  & 92.21\% \\
 & $10^{-4}$ 
 & 5.93  & 60.47\% \\
 & $10^{-5}$ 
 & 12.35 & 17.64\% \\
 & $10^{-6}$ 
 & 14.63 & 2.48\% \\
\hline
\end{tabular}
\caption{Application 1 - 1D concentrated case: relative error (eq. \eqref{eq:err2}) for all the configurations tested.}
\label{table:summary_Il_beta1}
\end{table}

\begin{figure}[htb!]
    \centering
    \begin{subfigure}[b]{0.48\textwidth}
        \centering
        \includegraphics[scale=0.48]{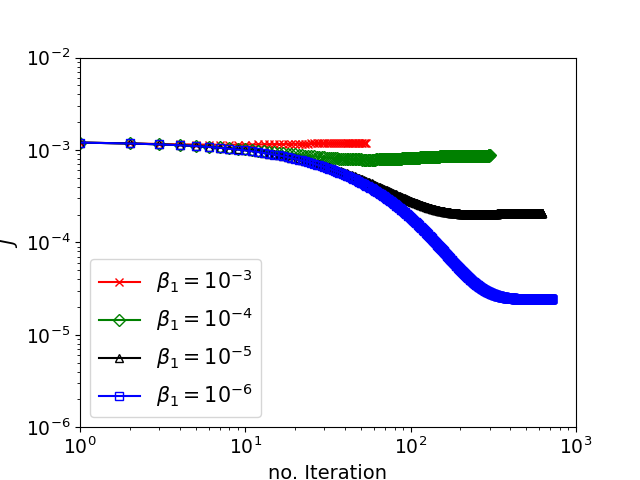}
    \end{subfigure}
    \hfill
    \begin{subfigure}[b]{0.48\textwidth}
        \centering
        \includegraphics[scale=0.48]{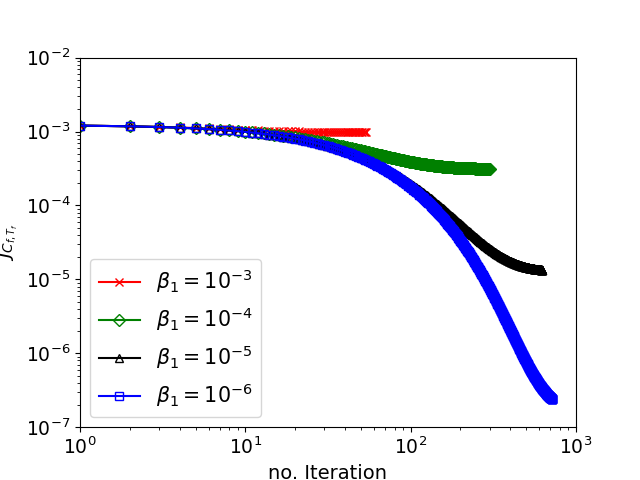}
    \end{subfigure}

    \begin{subfigure}[b]{0.48\textwidth}
        \centering
        \includegraphics[scale=0.48]{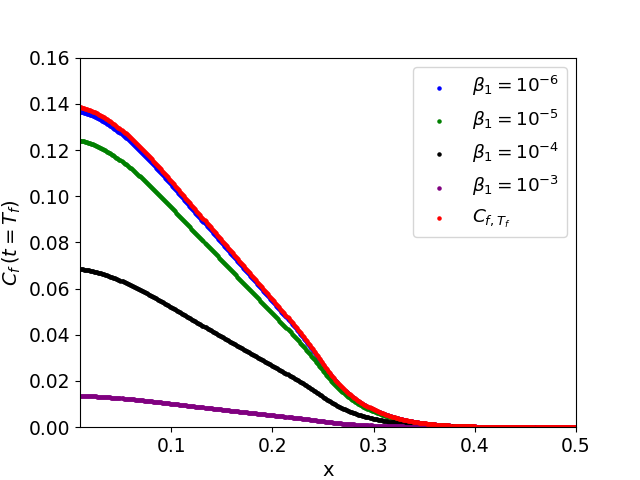}
    \end{subfigure}

     \caption{Application 1 - 1D concentrated case: convergence history of the objective function $J$ (top left), of the term $J_{c_{f,T_f}}$ (top right) and of the free drug concentration at the final time of the simulation (bottom) for different values of $\beta_1$ and for $I^0=5$.}
    \label{fig:plots_light_1D_I05_concentrated}
\end{figure}




\begin{figure}[htb!]
    \centering
    \begin{subfigure}[b]{0.48\textwidth}
        \centering
        \includegraphics[scale=0.48]{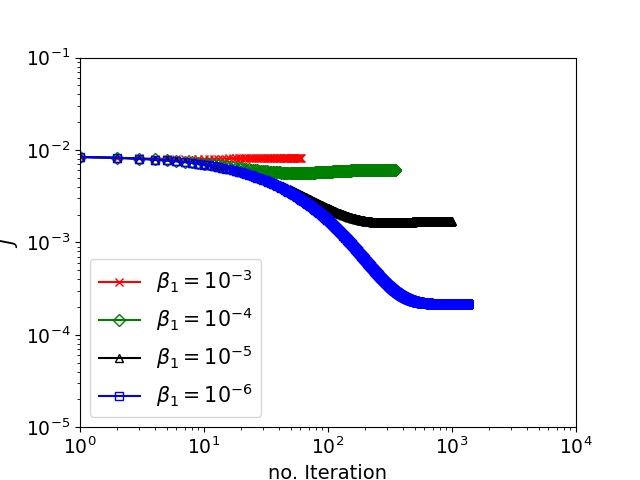}
    \end{subfigure}
    \hfill
    \begin{subfigure}[b]{0.48\textwidth}
        \centering
        \includegraphics[scale=0.48]{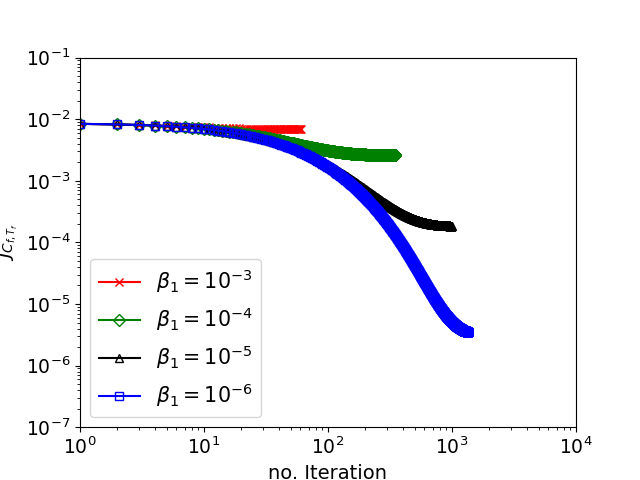}
    \end{subfigure}

    \begin{subfigure}[b]{0.48\textwidth}
        \centering
        \includegraphics[scale=0.48]{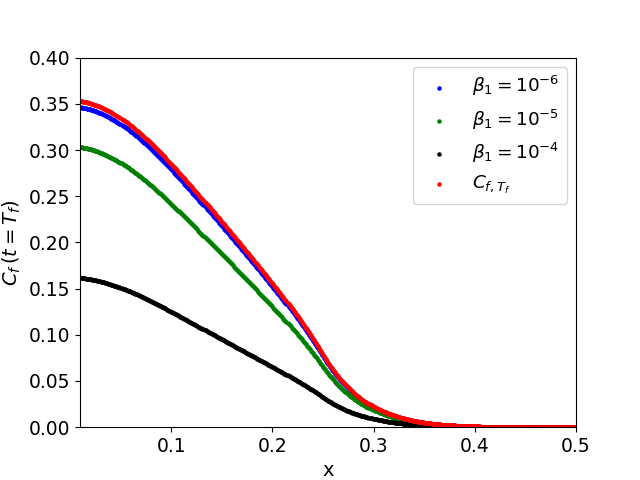}
    \end{subfigure}

    \caption{Application 1 - 1D concentrated case: convergence history of the objective function $J$ (top left), of the term $J_{c_{f,T_f}}$ (top right) and of the free drug concentration at the final time of the simulation for different values of $\beta_1$ (bottom) and for $I^0=15$.}
    \label{fig:plots_light_1D_I15_concentrated}
\end{figure}



\begin{figure}[htb!]
    \centering
    \begin{subfigure}[b]{0.48\textwidth}
        \centering
        \includegraphics[scale=0.48]{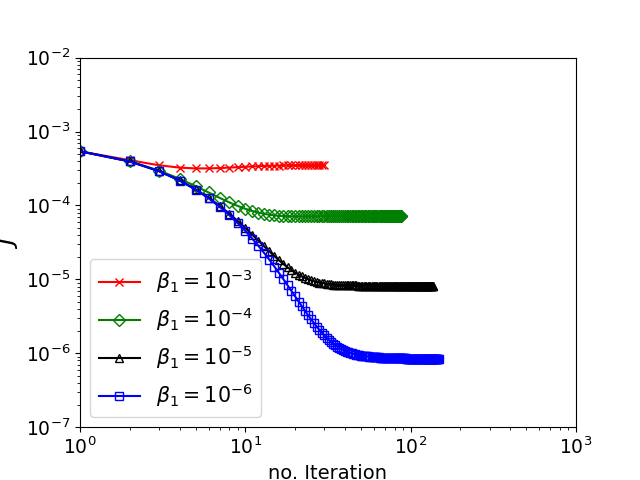}
    \end{subfigure}
    \hfill
    \begin{subfigure}[b]{0.48\textwidth}
        \centering
        \includegraphics[scale=0.48]{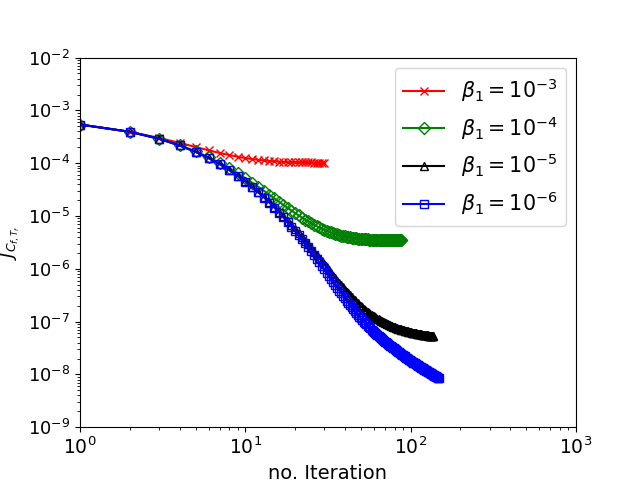}
    \end{subfigure}

    \begin{subfigure}[b]{0.48\textwidth}
        \centering
        \includegraphics[scale=0.48]{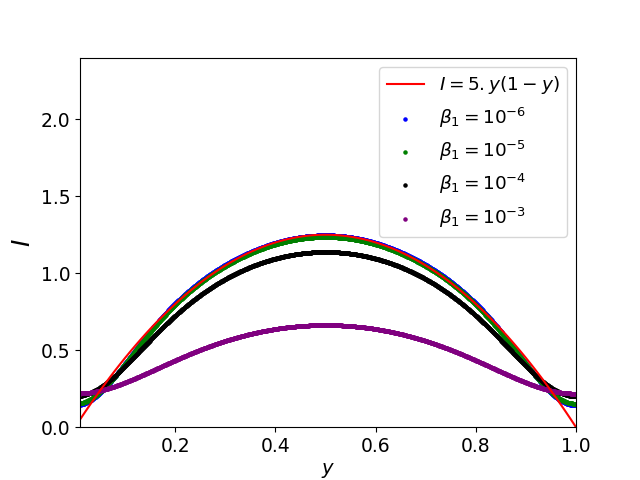}
    \end{subfigure}

    \caption{Application 1 - 2D concentrated case: convergence history of the objective function $J$ (top left), of the term $J_{c_{f,T_f}}$ and of the free drug concentration at the final time of the simulation for different values of $\beta_1$ (bottom) and for $I^0=5$.}
    \label{fig:plots_light_2D_I05_concentrated}
\end{figure}

\begin{figure}[htb!]
\centering
\resizebox{\textwidth}{!}{%
\begin{tabular}{cccccc}
     & $\beta_1=10^{-3}$ & $\beta_1=10^{-4}$ & $\beta_1=10^{-5}$ & $\beta_1=10^{-6}$ & \\
\raisebox{1.5cm}{$C_f(t=T_f)$} & \includegraphics[scale = 0.14]{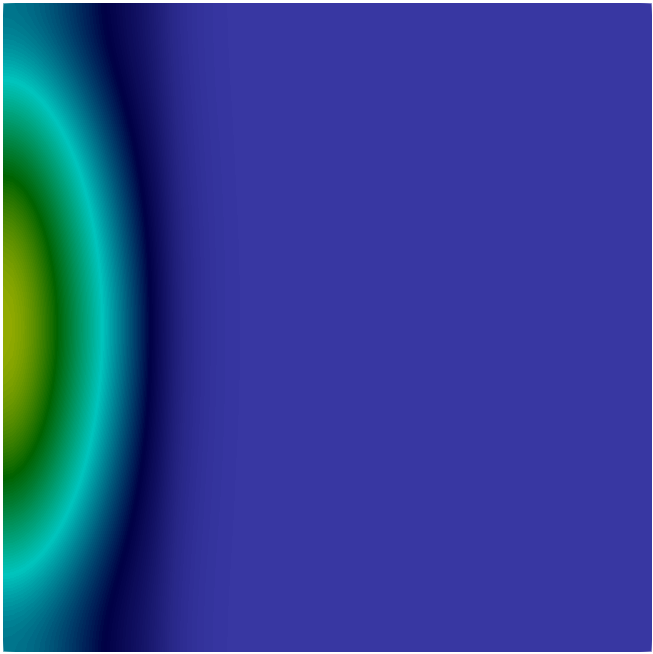} & \includegraphics[scale = 0.14]{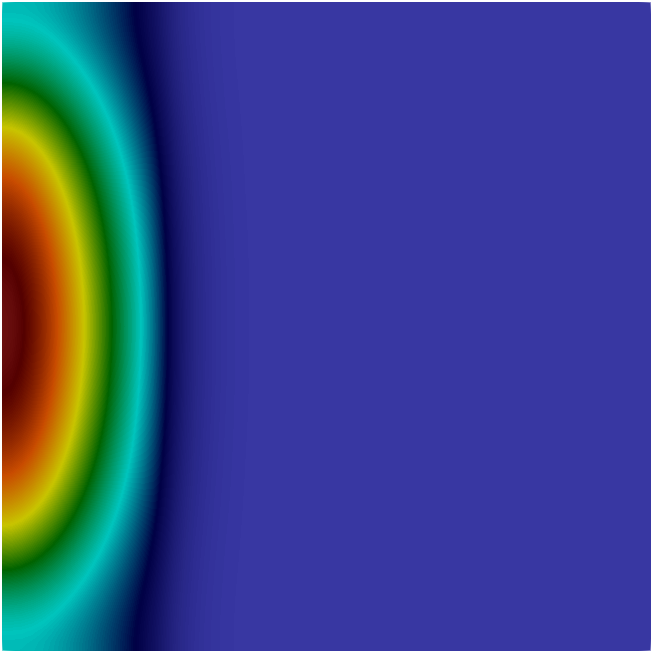} & \includegraphics[scale = 0.14]{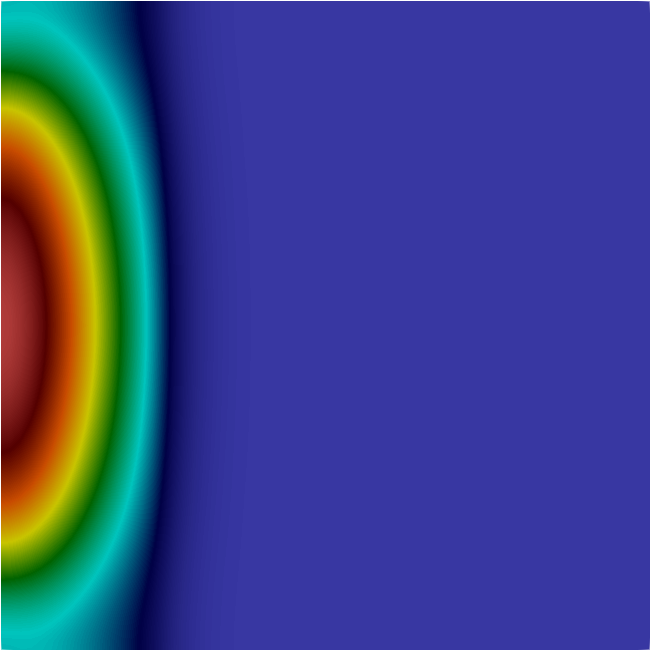} & \includegraphics[scale = 0.14]{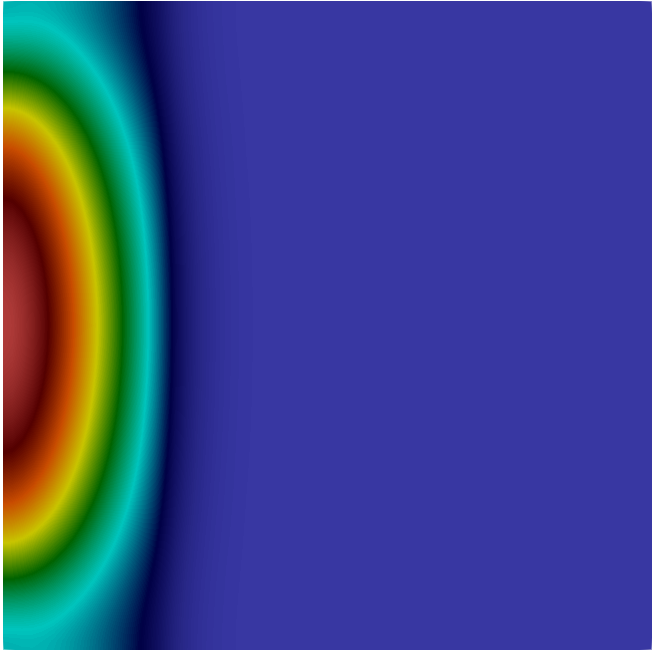} & \hspace{-0.4cm} \includegraphics[width=1.5cm, height=3.5cm]{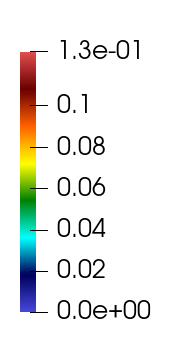} \\
\raisebox{1.5cm}{Error}  & \includegraphics[scale = 0.14]{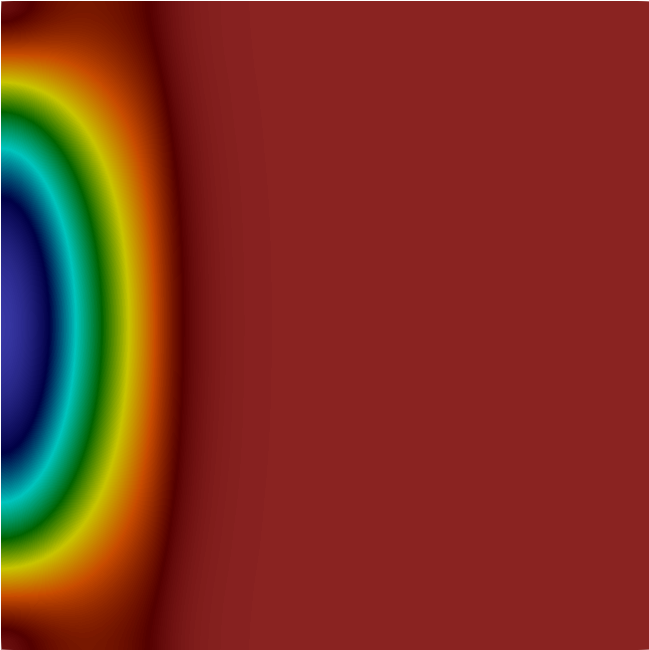} & \includegraphics[scale = 0.14]{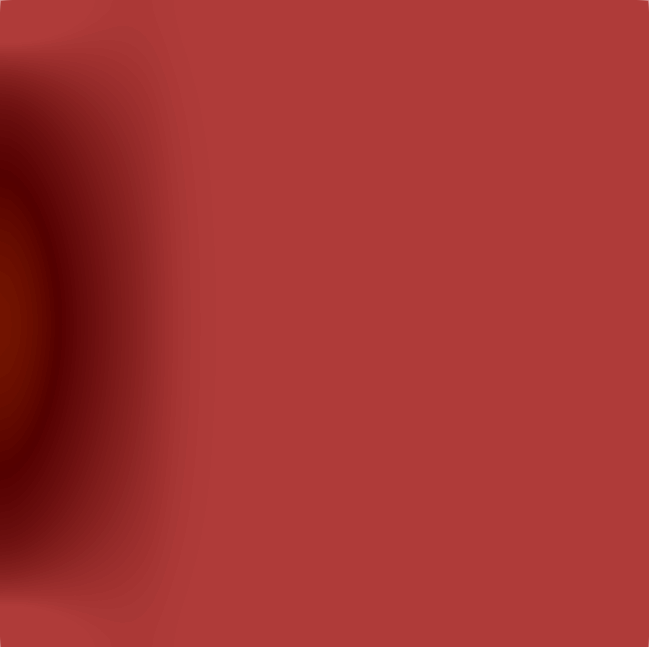} & \includegraphics[scale = 0.175]{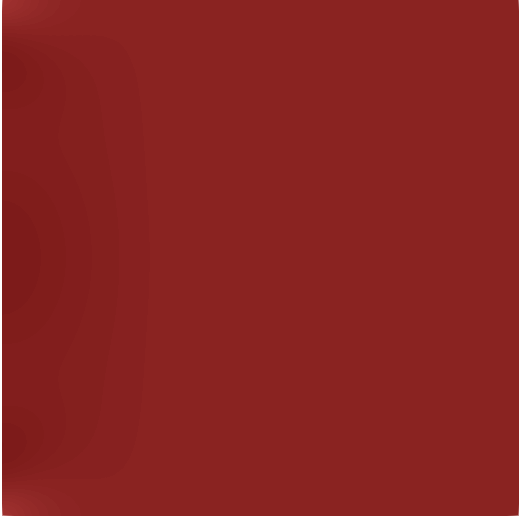} & \includegraphics[scale = 0.175]{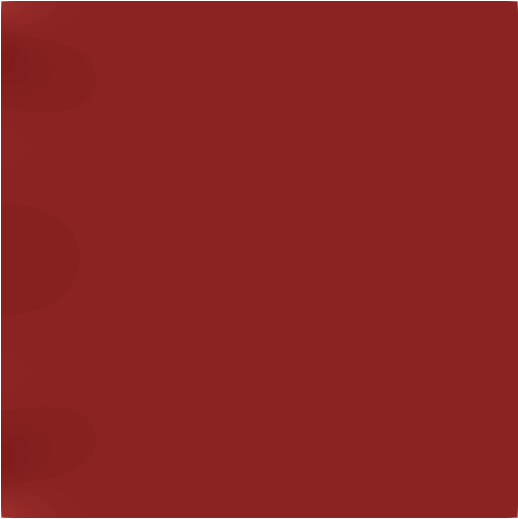} &  \hspace{-0.4cm} \includegraphics[width=1.5cm, height=3.5cm]{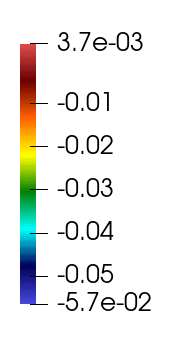}\\
\raisebox{1.5cm}{$C_{f,T_f}$}  & \multicolumn{4}{c}{\includegraphics[scale = 0.14]{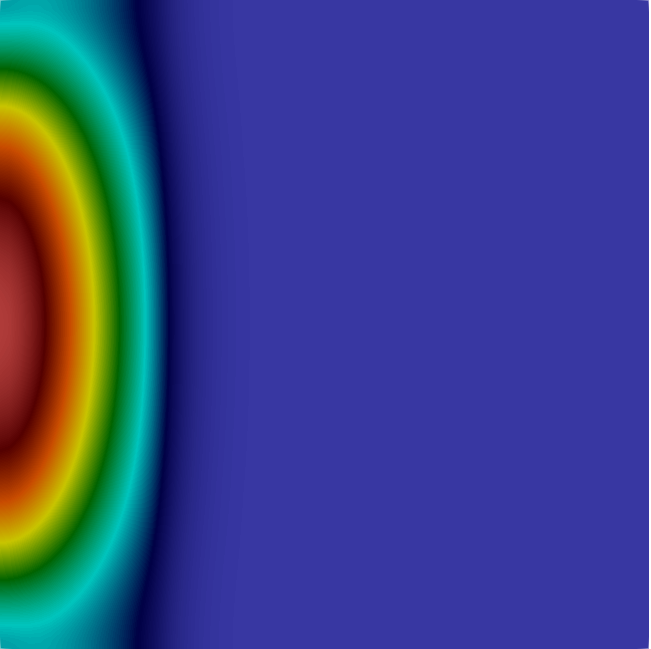}} & \hspace{-11.5cm} \includegraphics[width=1.5cm, height=3.5cm]{Cf_scale_I05.png} \\
\end{tabular}
}
\caption{Application 1 - 2D concentrated case: free drug
concentration at the final time of the simulation predicted by the OCP problem (first row), target distribution (third row) and  associated pointwise error (second row) for $I^0=5$.}
\label{fig:plots_light_2D_I05_concentrated_error}
\end{figure}

\begin{figure}[htb!]
    \centering
    \begin{subfigure}[b]{0.48\textwidth}
        \centering
        \includegraphics[scale=0.48]{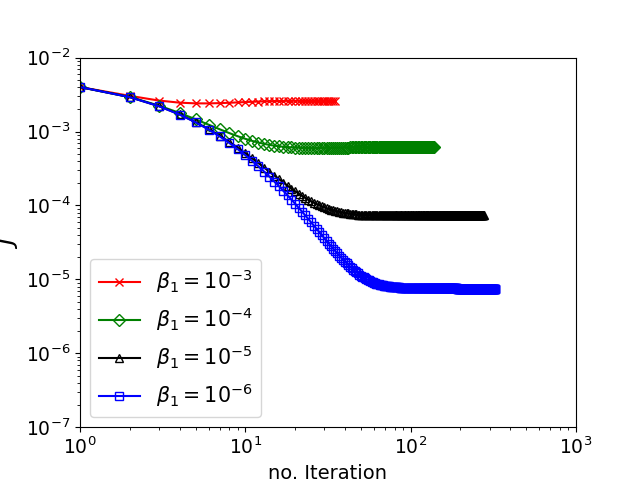}
    \end{subfigure}
    \hfill
    \begin{subfigure}[b]{0.48\textwidth}
        \centering
        \includegraphics[scale=0.48]{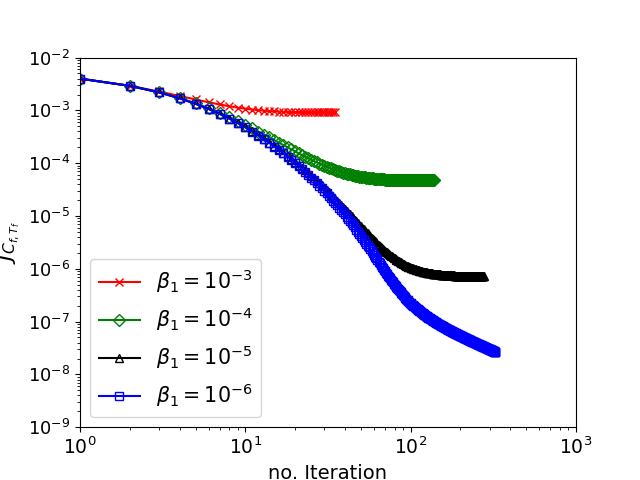}
    \end{subfigure}

    \begin{subfigure}[b]{0.48\textwidth}
        \centering
        \includegraphics[scale=0.48]{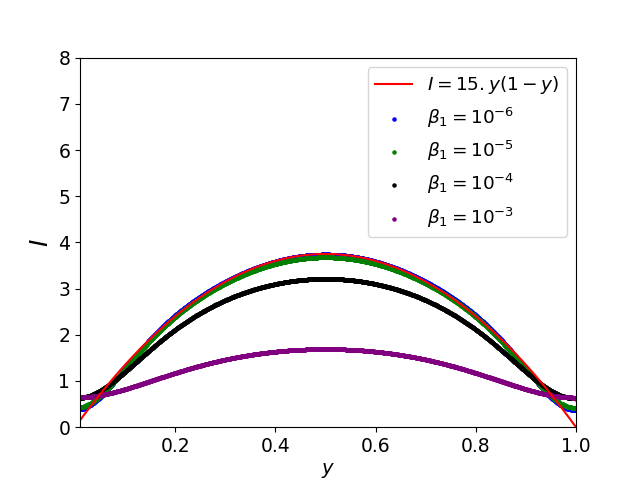}
    \end{subfigure}

    \caption{Application 1 - 2D concentrated case: convergence history of the objective function $J$ (top left), of the term $J_{c_{f,T_f}}$ (top right) and of the free drug concentration at the final time of the simulation (bottom) for different values of $\beta_1$ and for $I^0=15$.}
    \label{fig:plots_light_2D_I015_concentrated}
\end{figure}

\begin{figure}[htb!]
\centering
\resizebox{\textwidth}{!}{\begin{tabular}{cccccc}
     & $\beta_1=10^{-3}$ & $\beta_1=10^{-4}$ & $\beta_1=10^{-5}$ & $\beta_1=10^{-6}$ & \\
\raisebox{1.5cm}{$C_f(t=T_f)$} & \includegraphics[scale = 0.145]{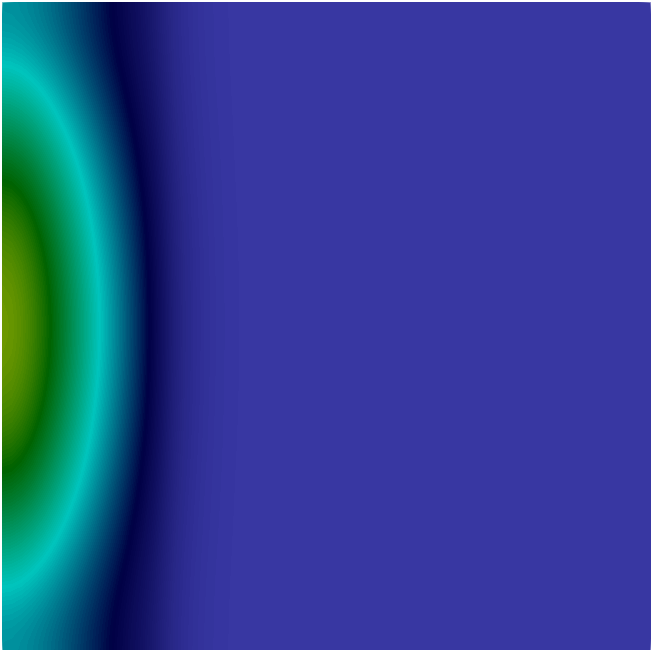} & \includegraphics[scale = 0.145]{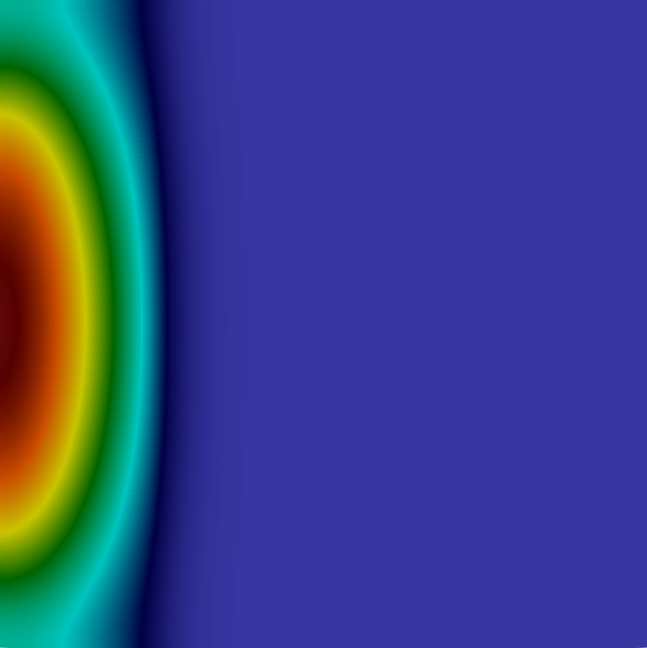} & \includegraphics[scale = 0.145]{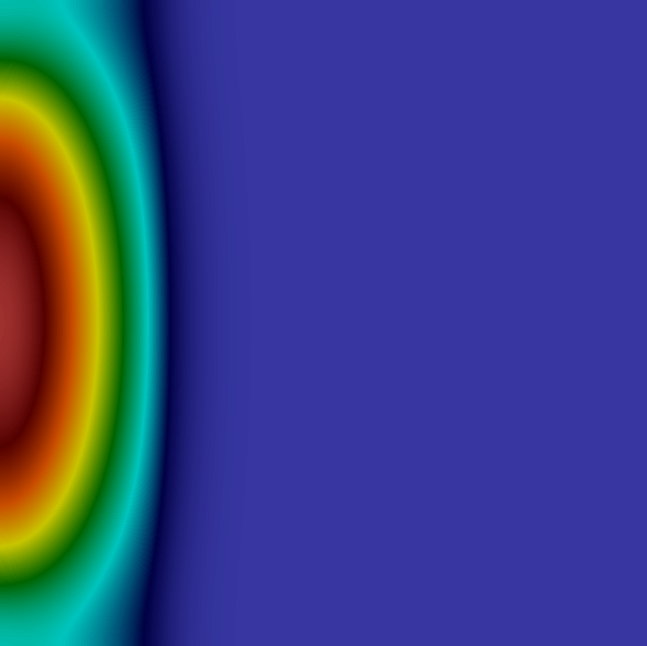} & \includegraphics[scale = 0.145]{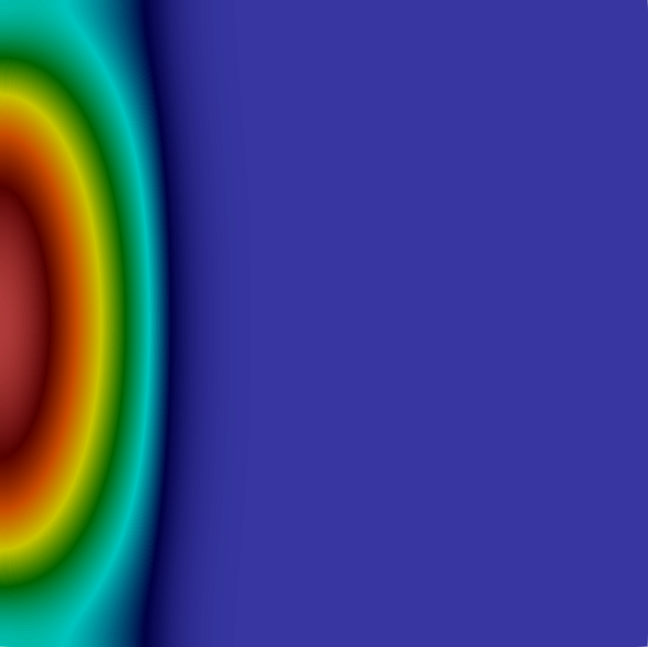} & \hspace{-0.7cm} \includegraphics[width=1.8cm, height=3.6cm]{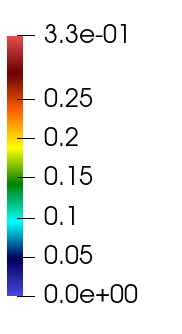} \\
\raisebox{1.5cm}{Error}  & \includegraphics[scale = 0.15]{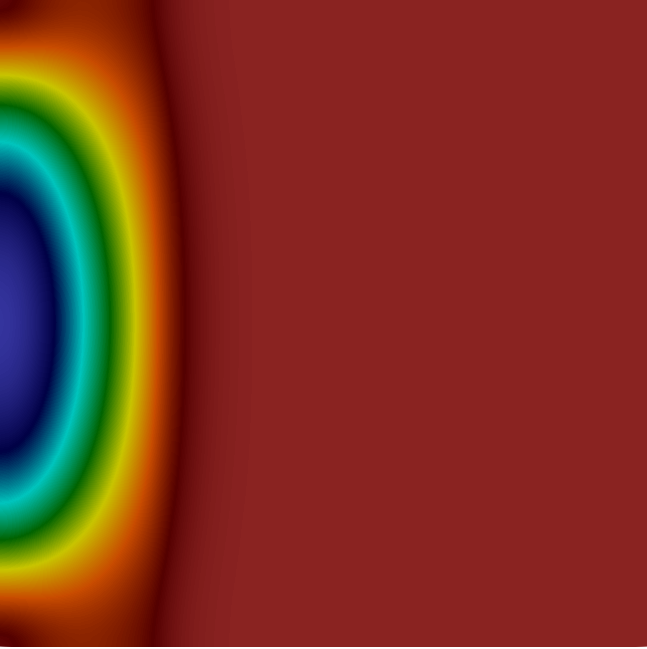} & \includegraphics[scale = 0.145]{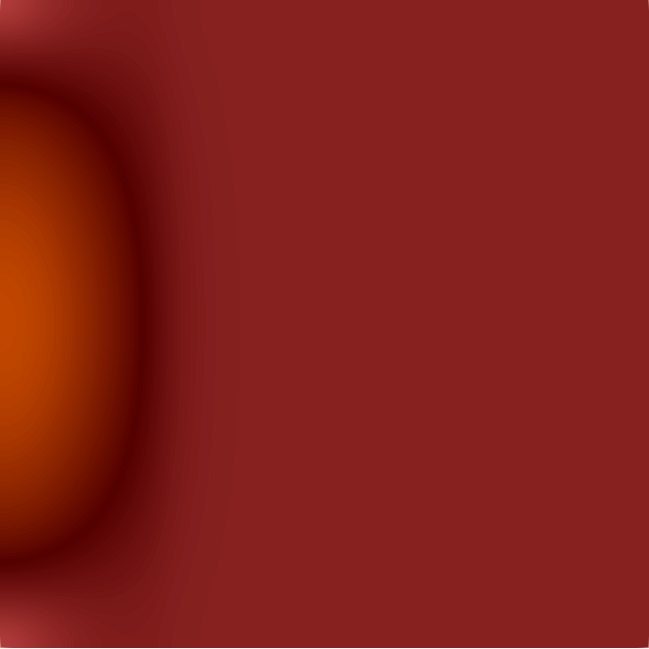} & \includegraphics[scale = 0.145]{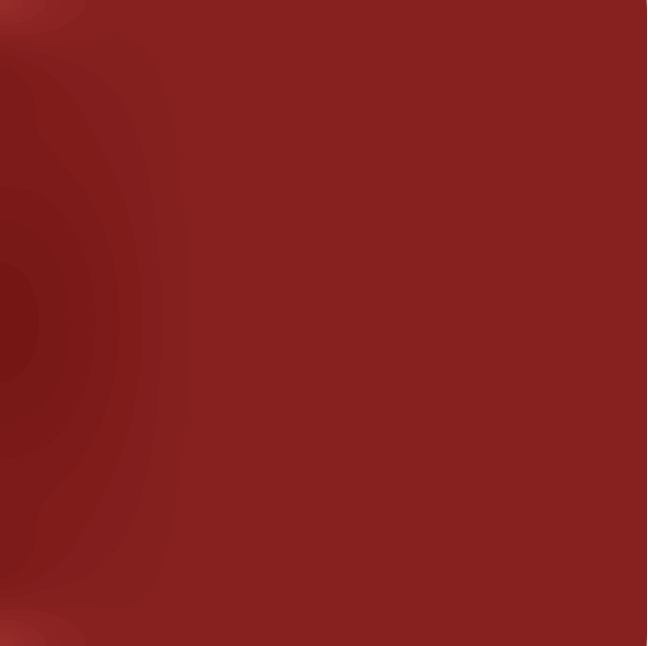} & \includegraphics[scale = 0.145]{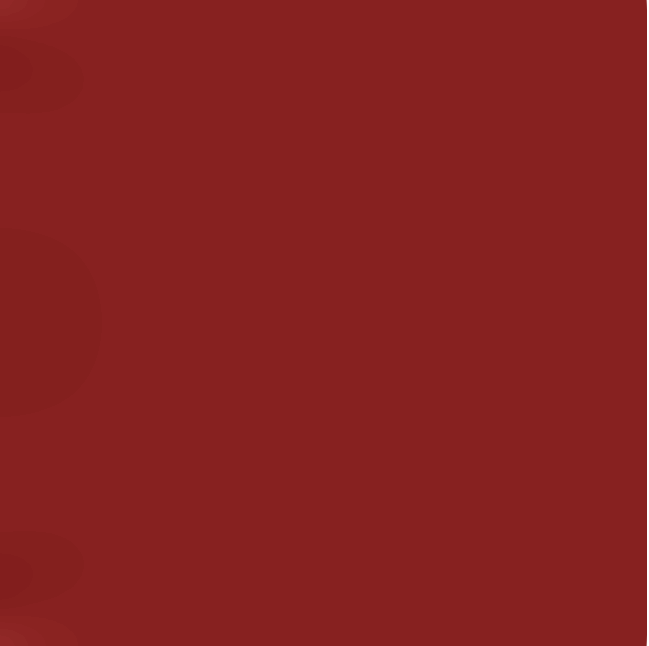} &  \hspace{-0.7cm} \includegraphics[width=1.8cm, height=3.7cm]{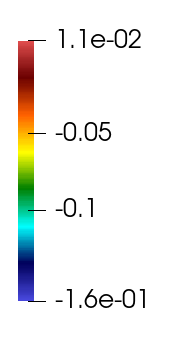}\\
\raisebox{1.5cm}{$C_{f,T_f}$}  & \multicolumn{4}{c}{\includegraphics[scale = 0.15]{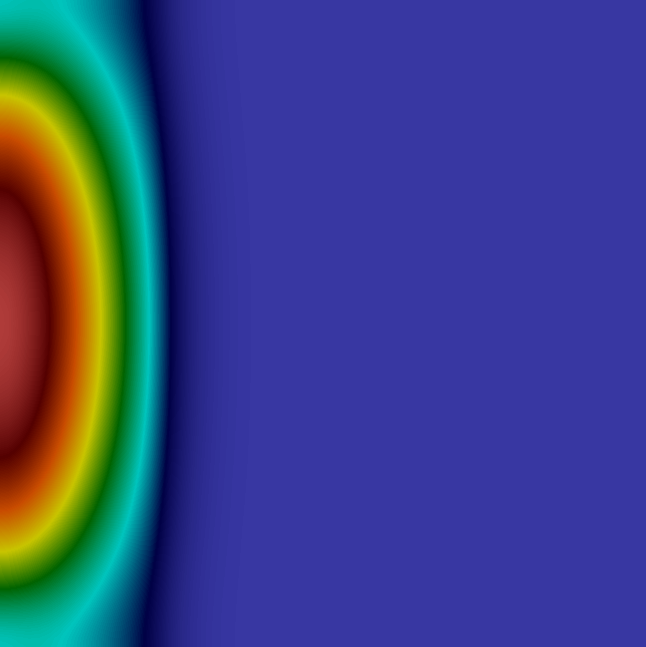}} & \hspace{-11.9cm} \includegraphics[width=1.8cm, height=3.6cm]{Cf_scale_I015_concentrated_2D.png} \\
\end{tabular}
}
\caption{Application 1 - 2D concentrated case: free drug
concentration at the final time of the simulation predicted by the OCP problem (first row), target distribution (third row) and associated pointwise error (second row) for $I^0=15$.}
\label{fig:plots_light_2D_I015_concentrated_error}
\end{figure}

\subsection{Application 2: Drug delivery in a coronary artery} \label{sec:drugcoronary}

The computational domain for this study is illustrated
in Figure \ref{fig:geometry}. 
The blood is described as a fluid with kinematic viscosity $\nu = 3.8 \cdot 10^{-6}$ m$^2$/s. In addition, the drug diffusivity $\varepsilon$ is set to $10^{-3}$. The mesh consists of 33469 cells. 
A periodic velocity waveform $V_i$ is applied at the inlet with a period of
$T_f = 0.8$ s (see Figure \ref{fig:geometry}). We set $\Delta t = 0.01$ s.  The catheter is assumed to inject the drug into the blood flow at the speed of $V_d = 4\cdot 10^{-2}$ m/s. 
It should be noted that, prior to the release of the drug, a fully periodic flow solution was attained. The drug is then released in the beginning of a period and is injected at a constant velocity thereafter.
Concerning the boundary and initial conditions for the drug concentration, we set $y_0 = y_\text{in} = y_{c} = y_\text{N} = 0$. 

Once we solved the Navier-Stokes equations system \eqref{tN-Steady} providing the flow field velocity $\mathbf{V}$, the validation is performed by using the same approach followed for the previous test cases, i.e., we generate data by solving the state problem \eqref{eq:general_transport} with an assigned value of $u$. We choose $u = 1$. Then we consider as target data $y_{T_f}$ the solution $y(t=T_f)$ at the final
time of the state problem simulation and we solve the OCP problem 
\eqref{eq:general_transport}, \eqref{eq:adjoint-cardio}, and \eqref{eq:optima-cardio}. In this case, no sensitivity analysis is performed with respect to $\beta_1$; instead,$\beta_1$ is directly set to $10^{-6}$, as our previous analysis proved this value to yield the best results.

The convergence features are reported in Table \ref{tab_cardio_3d}.  
Figure \ref{fig:plots_cardio_app_1} displays the spatial distribution of the drug at the final simulation time predicted by the OCP problem together with the
target distribution. We obtain a very good agreement from a qualitative viewpoint. Figure \ref{fig:plots_cardio_app_2} shows a close-up of the region next to the catheter where the drug variation is more pronounced. Here we also report the pointwise error which does not exceed $3.2 \cdot 10^{-3}$ by demonstrating a good accuracy of the OCP prediction. The error for the control variable is $E_{y} = |u - 1| = |1 - 0.9908| = 0.92\%$

\begin{table}[h!]
\centering
\begin{tabular}{|c|c|c|}
\hline
$N$ & \textbf{$J$} & \textbf{$J_{y_{T_f}}$} \\
\hline
10 & $3.2039 \cdot 10^{-12}$ & $6.0348 \cdot 10^{-13}$ \\
\hline
\end{tabular}
\caption{Application 2: convergence features.}\label{tab_cardio_3d}
\end{table}


\begin{figure}[htb!]
\centering
\begin{tabular}{cc}
\includegraphics[scale = 0.45]{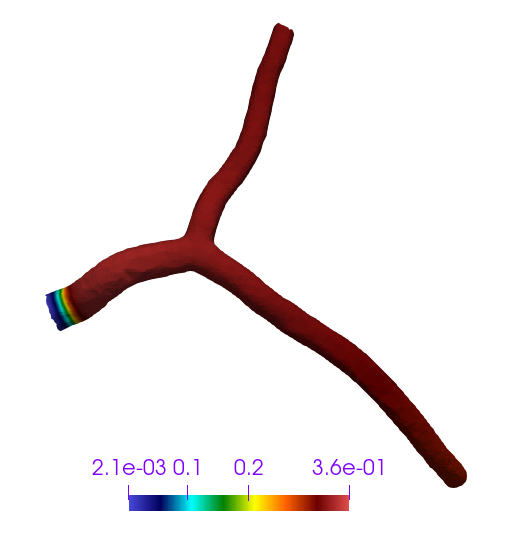} & \includegraphics[scale = 0.45]{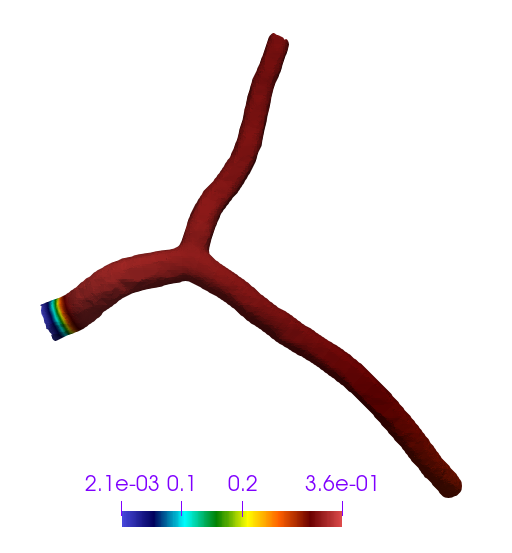} \\
\end{tabular}
\caption{Application 2: Drug concentration at the final time
of the simulation predicted by the OCP problem (left) and target distribution (right).}
\label{fig:plots_cardio_app_1}
\end{figure}

\begin{figure}[htb!]
\centering
\begin{tabular}{ccc}
\includegraphics[scale = 0.4]{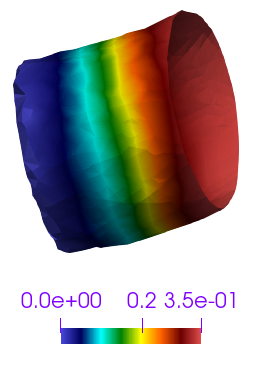} & \includegraphics[scale = 0.4]{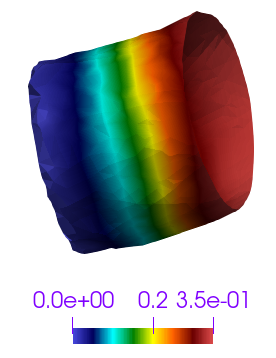} & 
\includegraphics[scale = 0.4]{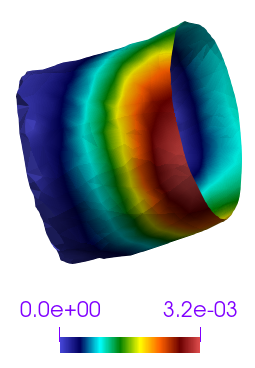}  \\
\end{tabular}
\caption{Application 2: Close-up of Figure \ref{fig:plots_cardio_app_1} (left and center) and point-wise error (right).}
\label{fig:plots_cardio_app_2}
\end{figure}





\section{Conclusions}\label{sec:concl}

In this work, we proposed a numerical approach based on the adjoint method to address optimal control problems for biomedical applications. After a preliminary validation against a well known academic benchmark \cite{fu2009}, our interest targeted two drug delivery problems, the former related to a light-driven system for cancer therapy (Application 1), the latter related to a catheter embedded in a patient-specific coronary artery (Application 2). 
Although both of these problems were already addressed in \cite{Ferreiral2022,ferreira2023prepub,ferreira2023, bazilevs2007} as  direct problems, here we reconsider them in an optimal control perspective in order to minimize the mismatch between numerical solution and target data through the tuning of a proper control variable. 

Numerical experiments not only demonstrated the accuracy of our approach, but also its flexibility and robustness in handling complex geometries (see Application 2), heterogeneous parameters (see the academic benchmark), and realistic boundary conditions (see both Application 1 and Application 2), highlighting its potential for the optimal design and control of complex biomedical systems. The computational outcomes have shown that lower values of the stabilization term weight ensure a higher performance in terms of accuracy.  
Anyway, we highlight a drawback in terms of efficiency. We have noted that the optimal control problem takes a computational time which is significantly higher than the direct one, about 3 to 1. This explains why numerical investigations of transient optimal control and inverse problems remain relatively limited. However,  we believe that by developing an appropriate Reduced Order Model (ROM) (see, e.g., \cite{ballarin2016, ballarin2017numerical, ballarinsupremizer2015}) the computational cost could be drastically reduced.  A natural extension of this work would be the development of a computational pipeline to tackle more complex optimal control problems in biomedical engineering, such as those arising in hemodynamics. Some relevant efforts in this context were recently introduced in \cite{rathore2025projection} but  much still remains to be done. 


\vskip 6mm
\noindent{\bf Acknowledgments}
\textbf{MG}, \textbf{ZM} and \textbf{GR} acknowledge the PRIN 2022 Project ``Machine learning for
fluid-structure interaction in cardiovascular problems: efficient solutions, model reduction, inverse problem (CUP
G53D23006830001)''. \textbf{PCA} and \textbf{GR} acknowledge the support provided by the European Union - NextGenerationEU, in the framework of the iNEST - Interconnected Nord-Est Innovation Ecosystem (iNEST ECS00000043 – CUP G93C22000610007) consortium and its CC5 Young Researchers initiative. The authors would also like to acknowledge INdAM-GNCS for its support.



\bibliographystyle{abbrv}
\bibliography{biblio.bib}

\begin{thebibliography}{10}

\bibitem{openFoam}
{O}pen{FOAM} {L}ibrary.
\newblock \url{https://openfoam.org/}.

\bibitem{ballarin2016}
F.~Ballarin, E.~Faggiano, S.~Ippolito, A.~Manzoni, A.~Quarteroni, G.~Rozza, and
  R.~Scrofani.
\newblock Fast simulations of patient-specific haemodynamics of coronary artery
  bypass grafts based on a {POD}-{G}alerkin method and a vascular shape
  parametrization.
\newblock {\em Journal of Computational Physics}, 315:609--628, 2016.

\bibitem{ballarin2017numerical}
F.~Ballarin, E.~Faggiano, A.~Manzoni, A.~Quarteroni, G.~Rozza, S.~Ippolito,
  C.~Antona, and R.~Scrofani.
\newblock {Numerical modeling of hemodynamics scenarios of patient-specific
  coronary artery bypass grafts}.
\newblock {\em Biomechanics and Modeling in Mechanobiology}, 16(4):1373--1399,
  2017.

\bibitem{ballarinsupremizer2015}
F.~Ballarin, A.~Manzoni, A.~Quarteroni, and G.~Rozza.
\newblock Supremizer stabilization of {POD}-{G}alerkin approximation of
  parametrized steady incompressible {N}avier-{S}tokes equations.
\newblock {\em International Journal for Numerical Methods in Engineering},
  102(5):1136--1161, 2015.

\bibitem{balzotti2024reduced}
C.~Balzotti, P.~Siena, M.~Girfoglio, G.~Stabile, J.~Due{\~n}as-Pamplona,
  J.~Sierra-Pallares, I.~Amat-Santos, and G.~Rozza.
\newblock {A reduced order model formulation for left atrium flow: an atrial
  fibrillation case}.
\newblock {\em Biomechanics and Modeling in Mechanobiology}, 23(4):1411--1429,
  2024.

\bibitem{barbeiro2017}
S.~Barbeiro, S.~Bardeji, J.~Ferreira, and L.~Pinto.
\newblock Non-fickian convection–diffusion models in porous media.
\newblock {\em Numerische Mathematik}, 2017.

\bibitem{bazilevs2007}
Y.~Bazilevs, V.~M. Calo, T.~E. Tezduyar, and T.~J.~R. Hughes.
\newblock Yz$\beta$ discontinuity capturing for advection‑dominated processes
  with application to arterial drug delivery.
\newblock {\em International Journal for Numerical Methods in Fluids},
  54(6–8):593--608, 2007.

\bibitem{bewley2001}
T.~R. Bewley.
\newblock Flow control: new challenges for a new renaissance.
\newblock {\em Progress in Aerospace Sciences}, 37(1):21--58, 2001.

\bibitem{dedeoptimal2007}
L.~Ded{\'e}.
\newblock Optimal flow control for {N}avier--{S}tokes equations: drag
  minimization.
\newblock {\em International Journal for Numerical Methods in Fluids},
  55(4):347--366, 2007.

\bibitem{d2012variational}
M.~D’Elia, M.~Perego, and A.~Veneziani.
\newblock {A variational data assimilation procedure for the incompressible
  Navier-Stokes equations in hemodynamics}.
\newblock {\em Journal of Scientific Computing}, 52(2):340--359, 2012.

\bibitem{Ferreiral2022}
J.~Ferreira, H.~Gomez, and L.~Pinto.
\newblock A numerical scheme for a model of drug delivery enhanced by light.
\newblock {\em International Conference on Mathematical Analysis and
  Applications in Science and Engineering ICMA2SC22}, 2022.

\bibitem{ferreira2023prepub}
J.~Ferreira, H.~Gómez, and L.~Pinto.
\newblock Finite differences–finite elements analysis and numerical
  simulation of a light-triggered drug delivery model.
\newblock {\em Pr\'e-Publica\c{c}\~{o}es do Departamento de Matem\'atica,
  Universidade de Coimbra, Preprint Number 21–15}.

\bibitem{ferreira2023}
J.~Ferreira, H.~Gómez, and L.~Pinto.
\newblock A numerical scheme for a partial romaro system motivated by
  light-triggered drug delivery.
\newblock {\em Applied Numerical Mathematics}, 184:101--120, 2023.

\bibitem{fevola2021}
E.~Fevola, F.~Ballarin, L.~Jim{\'e}nez-Juan, S.~Fremes, S.~Grivet-Talocia,
  G.~Rozza, and P.~Triverio.
\newblock An optimal control approach to determine resistance-type boundary
  conditions from in-vivo data for cardiovascular simulations.
\newblock {\em International Journal for Numerical Methods in Biomedical
  Engineering}, 37(10):e3516, 2021.

\bibitem{Formaggia}
L.~Formaggia, A.~Quarteroni, and A.~Veneziani.
\newblock {\em {Cardiovascular Mathematics: Modeling and simulation of the
  circulatory system}}, volume~1.
\newblock Springer Science \& Business Media, 2010.

\bibitem{fresca2020deep}
S.~Fresca, A.~Manzoni, L.~Ded{\`e}, and A.~Quarteroni.
\newblock {Deep learning-based reduced order models in cardiac
  electrophysiology}.
\newblock {\em PloS one}, 15(10):e0239416, 2020.

\bibitem{fu2009}
H.~Fu and H.~Rui.
\newblock A priori error estimates for optimal control problems governed by
  transient advection-diffusion equations.
\newblock {\em Journal of Scientific Computing}, 38:290--315, 2009.

\bibitem{Gad-elHak2003}
M.~Gad-el Hak, A.~Pollard, and J.-P. Bonnet.
\newblock {\em Flow control: Fundamentals and practices}, volume~53.
\newblock Springer Science \& Business Media, 2003.

\bibitem{girfogliononintrusive}
M.~Girfoglio, L.~Scandurra, F.~Ballarin, G.~Infantino, F.~Nicolo, A.~Montalto,
  G.~Rozza, R.~Scrofani, M.~Comisso, and F.~Musumeci.
\newblock Non-intrusive data-driven {ROM} framework for hemodynamics problems.
\newblock {\em Acta mechanica sinica}, 37(7):1183--1191, 2021.

\bibitem{gunzburger2002}
M.~D. Gunzburger.
\newblock {\em Perspectives in flow control and optimization}, volume~5 of {\em
  Advances in Design and Control}.
\newblock Society for Industrial and Applied Mathematics (SIAM), Philadelphia,
  PA, 2003.

\bibitem{hinze2009}
M.~Hinze, R.~Pinnau, M.~Ulbrich, and S.~Ulbrich.
\newblock {\em Optimization with PDE constraints}, volume~23.
\newblock Springer Science \& Business Media, 2008.

\bibitem{hounumerical1999}
L.~Hou and S.~Ravindran.
\newblock Numerical approximation of optimal flow control problems by a penalty
  method: Error estimates and numerical results.
\newblock {\em SIAM Journal on Scientific Computing}, 20(5):1753--1777, 1999.

\bibitem{ito2008}
K.~Ito and K.~Kunisch.
\newblock {\em Lagrange multiplier approach to variational problems and
  applications}.
\newblock SIAM, 2008.

\bibitem{ji2019}
Y.~Ji et~al.
\newblock A light-facilitated drug delivery system from a
  pseudo-protein/hyaluronic acid nanocomplex with improved anti-tumor effects.
\newblock {\em Nanoscale}, 2019.

\bibitem{kalaydina2018}
R.~V. Kalaydina et~al.
\newblock Recent advances in “smart” delivery systems for extended drug
  release in cancer therapy.
\newblock {\em International Journal of Nanomedicine}, 2018.

\bibitem{kelley1999iterative}
C.~T. Kelley.
\newblock {\em Iterative methods for optimization}.
\newblock SIAM, 1999.

\bibitem{lassila2013reduced}
T.~Lassila, A.~Manzoni, A.~Quarteroni, and G.~Rozza.
\newblock {A reduced computational and geometrical framework for inverse
  problems in hemodynamics}.
\newblock {\em International Journal for Numerical Methods in Biomedical
  Engineering}, 29(7):741--776, 2013.

\bibitem{manzoni2012model}
A.~Manzoni, A.~Quarteroni, and G.~Rozza.
\newblock {Model reduction techniques for fast blood flow simulation in
  parametrized geometries}.
\newblock {\em International Journal for Numerical Methods in Biomedical
  Engineering}, 28(6-7):604--625, 2012.

\bibitem{quarteroni2009}
A.~Quarteroni and S.~Quarteroni.
\newblock {\em Numerical models for differential problems}, volume~2.
\newblock Springer, 2009.

\bibitem{quarteronioptimal}
A.~Quarteroni and G.~Rozza.
\newblock Optimal control and shape optimization of aorto-coronaric bypass
  anastomoses.
\newblock {\em Mathematical Models and Methods in Applied Sciences},
  13(12):1801--1823, 2003.

\bibitem{quarteroni2016geometric}
A.~Quarteroni, A.~Veneziani, and C.~Vergara.
\newblock {Geometric multiscale modeling of the cardiovascular system, between
  theory and practice}.
\newblock {\em Computer Methods in Applied Mechanics and Engineering},
  302:193--252, 2016.

\bibitem{rathore2025projection}
S.~Rathore, P.~C. Africa, F.~Ballarin, F.~Pichi, M.~Girfoglio, and G.~Rozza.
\newblock {Projection-based reduced order modelling for unsteady parametrized
  optimal control problems in 3D cardiovascular flows}.
\newblock {\em Computer Methods and Programs in Biomedicine}, page 108813,
  2025.

\bibitem{Romarowski2018}
R.~M. Romarowski, A.~Lefieux, S.~Morganti, A.~Veneziani, and F.~Auricchio.
\newblock Patient-specific {CFD} modelling in the thoracic aorta with
  {PC-MRI}--based boundary conditions: A least-square three-element
  {W}indkessel approach.
\newblock {\em International journal for numerical methods in biomedical
  engineering}, 34(11):e3134, 2018.

\bibitem{sankaran2012patient}
S.~Sankaran, M.~Esmaily~Moghadam, A.~M. Kahn, E.~E. Tseng, J.~M. Guccione, and
  A.~L. Marsden.
\newblock {Patient-specific multiscale modeling of blood flow for coronary
  artery bypass graft surgery}.
\newblock {\em Annals of biomedical engineering}, 40(10):2228--2242, 2012.

\bibitem{siena2025hybrid}
P.~Siena, P.~C. Africa, M.~Girfoglio, and G.~Rozza.
\newblock A hybrid reduced order model to enforce outflow pressure boundary
  conditions in computational haemodynamics.
\newblock {\em Biomechanics and Modeling in Mechanobiology}, 2025.

\bibitem{sienadata}
P.~Siena, M.~Girfoglio, F.~Ballarin, and G.~Rozza.
\newblock {Data-driven reduced order modelling for patient-specific
  hemodynamics of coronary artery bypass grafts with physical and geometrical
  parameters}.
\newblock {\em Journal of Scientific Computing}, 94(2):38, 2023.

\bibitem{strazzullo2021model}
M.~Strazzullo.
\newblock {\em {Model Order Reduction for Nonlinear and Time-Dependent
  Parametric Optimal Flow Control Problems}}.
\newblock Master thesis. SISSA, 2021.

\bibitem{tang2019}
M.~Tang et~al.
\newblock Can intracellular drug delivery using hyaluronic acid functionalised
  ph-sensitive liposomes overcome gemcitabine resistance in pancreatic cancer?
\newblock {\em Journal of Controlled Release}, 2019.

\bibitem{thomas2021patient}
B.~Thomas, K.~Sumam, and N.~Sajikumar.
\newblock {Patient specific modelling of blood flow in coronary artery}.
\newblock {\em Journal of Applied Fluid Mechanics}, 14(5):1469--1482, 2021.

\bibitem{viceconti2015big}
M.~Viceconti, P.~Hunter, and R.~Hose.
\newblock {Big data, big knowledge: big data for personalized healthcare}.
\newblock {\em IEEE journal of biomedical and health informatics},
  19(4):1209--1215, 2015.

\bibitem{wang2020}
X.~Wang et~al.
\newblock Near-infrared photoresponsive drug delivery nanosystems for cancer
  photo-chemotherapy.
\newblock {\em Journal of Nanobiotechnology}, 2020.

\bibitem{zakia2021}
Z.~Zainib, F.~Ballarin, S.~Fremes, P.~Triverio, L.~Jim\'{e}nez-Juan, and
  G.~Rozza.
\newblock Reduced order methods for parametric optimal flow control in coronary
  bypass grafts, toward patient-specific data assimilation.
\newblock {\em International Journal for Numerical Methods in Biomedical
  Engineering}, 37(12):Paper No. e3367, 2021.

\end{thebibliography}

\end{document}